\pdfoutput=1
\documentclass[preliminary]{eptcs}
\providecommand{\event}{GCM 2021}%
\hypersetup{ colorlinks
           , urlcolor=RoyalBlue
           , linkcolor=Red
           , citecolor=Green
           }

\usepackage{graphicx}
\usepackage[dvipsnames]{xcolor}

  \usepackage{tikz}
  \usetikzlibrary{
    cd,
    petri,
    backgrounds,
    arrows,
    positioning,
    decorations.markings,
    calc,
    fit,
    shapes.multipart
  }

\usepackage{mathtools}
\usepackage{amssymb}

\usepackage{amsthm}
\usepackage{stmaryrd}

\usepackage{relsize}
\usepackage{microtype}

\usepackage{multicol}
\usepackage{csquotes}
\usepackage{xspace}

\usepackage{enumitem}
\usepackage{fontawesome}

\usepackage{tabularx, cellspace}
\usepackage[capitalize]{cleveref}

\usepackage{url}
  \newcounter{theoremUnified} 
  \numberwithin{theoremUnified}{section} 
  \numberwithin{theoremUnified}{section} 

  \newtheoremstyle{plainStyle} 
  {2mm} 
  {2mm} 
  {} 
  {} 
  {\bfseries} 
  {.} 
  {.5em} 
  {} 

  \newtheoremstyle{italicStyle} 
  {2mm} 
  {2mm} 
  {\itshape} 
  {} 
  {\bfseries} 
  {.} 
  {.5em} 
  {} 

\newcommand{\Naturals}{\mathbb{N}} 

\def\backgrnd{black!10}	

\pgfdeclarelayer{edgelayer}
\pgfdeclarelayer{nodelayer}
\pgfdeclarelayer{background}
\pgfsetlayers{background, edgelayer,nodelayer,main}

\tikzstyle{place}=
[circle,thick,draw=blue!75,fill=blue!20,minimum size=6mm]
\tikzstyle{antiplace}=
[circle,thick,draw=red!75,fill=red!20,minimum size=6mm]
\tikzstyle{transition}=
[rectangle,thick,draw=black!75,fill=black!20,minimum size=4mm]
\tikzstyle{inarrow}=[->, >=stealth, shorten >=.03cm,line width=1.5]
\tikzstyle{antiarrow}=[<-, red!75,  >=stealth, shorten >=.03cm,line width=1.5]
\tikzset{NodeBox/.style = {rectangle, draw=black!50, fill=black!10, rounded corners, minimum height=30, minimum width=30, align=center, font=\tiny\ttfamily}}
\tikzset{ArrowBox/.style = {rectangle, draw=black!50, fill=black!10, minimum width=40, align=left, font=\tiny\ttfamily}}

\tikzset{
	pics/netA/.style args={#1/#2/#3/#4/#5/#6/#7}{code={

					\node [place,label=above:$p_1$, tokens={%
								#1
							}] (-pl_1) {};

					\node [transition,label=above:$t$, label=below:#5] (-tr_1) [right = of -pl_1] {};

					\node [place,label=above:$p_2$, tokens={%
								#2
							}] (-pl_2) [right = of -tr_1] {};

					\node [transition,label=left:$v$, label=above:#6] (-tr_2) [below = of -tr_1] {};
					\node [transition,label=below:$u$, label=above:#7] (-tr_3) [below = of -tr_2] {};

					\node [place,label=below:$p_3$, tokens={%
								#3
							}] (-pl_3) [left = of -tr_3] {};

					\node [place,label=below:$p_4$, tokens={%
								#4
							}] (-pl_4) [right = of -tr_3] {};

					\draw[inarrow] (-pl_1) -- (-tr_1);
					\draw[inarrow] (-tr_1) -- (-pl_2);
					\draw[inarrow] (-pl_2) -- (-tr_2);
					\draw[inarrow] (-tr_2) -- (-pl_3);
					\draw[inarrow] (-tr_2) -- (-pl_4);
					\draw[inarrow] (-pl_3) -- (-tr_3);
					\draw[inarrow] (-tr_3) -- (-pl_4);
				}}
}

\newcommand{\Msets}[1]{{#1}^{\oplus}} 

\newcommand{\cp}{\fatsemi} 

\newcommand{\Obj}[1]{\operatorname{Obj} \, #1} 

\newcommand{\Hom}[3]{\operatorname{Hom}_{\,#1}\left[#2,#3\right]} 

\newcommand{\CategoryC}{\mathcal{C}}

\newcommand{\Fun}[1]{{#1}^\sharp} 
\newcommand{\NetSem}[1]{\left( #1, \Fun{#1}\right)} 

\newcommand{\Semantics}{\mathcal{S}} 

\newcommand{\Free}[1]{\mathfrak{F}\left(#1\right)} 

\newcommand{\Grothendieck}[1]{\textstyle\int{#1}} 
\newcommand{\GrothendieckS}[1]{\Grothendieck{\Fun{#1}}} 

\newcommand{\harpvecsign}{\scriptscriptstyle\rightharpoonup}
\newcommand{\harpoonvec}[2]{%
  \ifx\displaystyle#1\doalign{$\harpvecsign$}{#1#2}\fi
  \ifx\textstyle#1\doalign{$\harpvecsign$}{#1#2}\fi
  \ifx\scriptstyle#1\doalign{\scalebox{.6}[.9]{$\harpvecsign$}}{#1#2}\fi
  \ifx\scriptscriptstyle#1\doalign{\scalebox{.5}[.8]{$\harpvecsign$}}{#1#2}\fi
}
\newcommand{\doalign}[2]{%
 {\vbox{\offinterlineskip\ialign{\hfil##\hfil\cr#1\cr$#2$\cr}}}%
}

\newcommand{\LiftSpan}[1]{\widehat{#1}} 

\let\cate\mathbf
\newcommand{\Set}{\cate{Set}} 
\newcommand{\Span}{\cate{Span}} 
\newcommand{\Cat}{\cate{Cat}} 
\newcommand{\Petri}{\cate{Petri}} 
\newcommand{\FSSMC}{\cate{FSSMC}} 
\def\fssmc{\textsc{fssmc}\xspace}
\newcommand{\PetriS}[1]{\cate{Petri}^{#1}} 
  \newcommand{\PetriSpan}{\PetriS{\Span}} 
  \newcommand{\PetriHier}{\PetriS{\tikz\fill[black]
  (0,0)  -- +(.1,0) -- +(60:.1) -- cycle
  (.1,0)  -- +(.1,0) -- +(60:.1) -- cycle
  (60:.1) -- +(.1,0) -- +(60:.1) -- cycle
;}} 
  \newcommand{\Mon}{\cate{Mon}} 

\newcommand{\Arrow}[1]{#1^{\to}}
\newcommand{\Object}[1]{#1^{\bullet}}
\newcommand{\Dom}[1]{\mathrm{dom}_{#1}}
\newcommand{\Cod}[1]{\mathrm{cod}_{#1}}

\newcommand{\Play}[1]{\blacktriangleright_{#1}}
\newcommand{\Stop}[1]{\blacksquare_{#1}}
\newcommand{\Suchthat}[2]{\left\{#1 \: \middle\vert \: #2\right\}} 

\tikzset{ 
  oriented WD/.style={
    every to/.style={
      out=0,in=180,draw
    },
    label/.style={
      font=\everymath\expandafter{\the\everymath\scriptstyle},
      inner sep=0pt,
      node distance=2pt and -2pt
    },
    semithick,
    node distance=1 and 1,
    decoration={
      markings, mark=at position \stringdecpos with \stringdec
    },
    ar/.style={
      postaction={decorate}
    },
    execute at begin picture={
      \tikzset{
        x=\bbx, y=\bby,
        every fit/.style={
          inner xsep=\bbx, inner ysep=\bby
        }
      }
    }
  },
  string decoration/.store in=\stringdec,
  string decoration={
    \arrow{stealth};
  },
  string decoration pos/.store in=\stringdecpos,
  string decoration pos=.7,
  bbx/.store in=\bbx,
  bbx = 1.5cm,
  bby/.store in=\bby,
  bby = 1.5ex,
  bb port sep/.store in=\bbportsep,
  bb port sep=1.5,
  bb port length/.store in=\bbportlen,
  bb port length=4pt,
  bb penetrate/.store in=\bbpenetrate,
  bb penetrate=0,
  bb min width/.store in=\bbminwidth,
  bb min width=1cm,
  bb rounded corners/.store in=\bbcorners,
  bb rounded corners=2pt,
  bb small/.style={
    bb port sep=1,
    bb port length=2.5pt,
    bbx=.4cm, bb min width=.4cm,
    bby=.7ex
  },
  bb medium/.style={
    bb port sep=1,
    bb port length=2.5pt,
    bbx=.4cm,
    bb min width=.4cm,
    bby=.9ex
  },
  bb/.code 2 args={
    \pgfmathsetlengthmacro{\bbheight}{\bbportsep * (max(#1,#2)+1) * \bby}
    \pgfkeysalso{
      draw,
      minimum height=\bbheight,
      minimum width=\bbminwidth,
      outer sep=0pt,
      rounded corners=\bbcorners,
      thick,
      prefix after command={
        \pgfextra{\let\fixname\tikzlastnode}
      },
      append after command={
        \pgfextra{
          \draw
          \ifnum #1=0
            {}
          \else
            foreach \i in {1,...,#1} {
              ($(\fixname.north west)!{\i/(#1+1)}!(\fixname.south west)$) +(-
            \bbportlen,0)
            coordinate (\fixname_in\i) -- +(\bbpenetrate,0) coordinate (\fixname_in\i')
            }
          \fi
          \ifnum
            #2=0{}
          \else
            foreach \i in {1,...,#2} {
            ($(\fixname.north east)!{\i/(#2+1)}!(\fixname.south east)$) +(-
            \bbpenetrate,0)
            coordinate (\fixname_out\i') -- +(\bbportlen,0) coordinate (\fixname_out\i)
            }
          \fi;
        }
      }
    }
  },
  bb name/.style={
    append after command={
      \pgfextra{
        \node[anchor=north] at (\fixname.north) {#1}
      ;}
    }
  }
}

\newcommand{\from}{\leftarrow}
\newcommand{\From}[1]{\xleftarrow{#1}}
\newcommand{\To}[1]{\xrightarrow{#1}}

\def\id{\mathrm{id}}
\usepackage[all, 2cell]{xy}\UseAllTwocells
\usepackage{adjustbox}
\newenvironment{adju}[1][0.925]{%
\begin{center}\begin{adjustbox}{max height=0.5\textheight, max width=#1\textwidth}}{\end{adjustbox}\end{center}}

\newtheoremstyle{exampstyle}
  {2mm}
  {2mm}
  {}
  {}
  {\bfseries}
  {.}
  {.5em}
  {}
\newtheorem{theorem}{Theorem}
 \AfterEndEnvironment{theorem}{\noindent\ignorespaces}

 \AfterEndEnvironment{proposition}{\noindent\ignorespaces}
\newtheorem{definition}{Definition}
  \AfterEndEnvironment{definition}{\noindent\ignorespaces}

 \AfterEndEnvironment{fact}{\noindent\ignorespaces}
\AfterEndEnvironment{proof}{\noindent\ignorespaces}
\theoremstyle{exampstyle}
\newtheorem{notation}{Notation}{\itshape}{\rmfamily}
  \AfterEndEnvironment{notation}{\noindent\ignorespaces}
{\itshape}{\rmfamily}
  \AfterEndEnvironment{counterexample}{\noindent\ignorespaces}
\newtheorem{remark}{Remark}
  \AfterEndEnvironment{remark}{\noindent\ignorespaces}
\newtheorem{example}{Example}
  \AfterEndEnvironment{example}{\noindent\ignorespaces}
\title{A Categorical Semantics for Hierarchical Petri Nets}
\author{Fabrizio Romano Genovese%
\thanks{The first author was supported by the project MIUR PRIN 2017FTXR7S “IT-MaTTerS” and by the \href{https://gitcoin.co/grants/1086/independent-ethvestigator-program}{Independent Ethvestigator Program}.}
  \email{0000-0001-7792-1375}%
  \institute{University of Pisa / Statebox}%
  \email{fabrizio.romano.genovese@gmail.com}
  \and
  Jelle Herold%
  \email{0000-0002-1966-2536}%
  \institute{Statebox}%
  \email{research@statebox.io}
  \and
  Fosco Loregian%
  \thanks{The third author was supported by the ESF funded Estonian IT Academy research measure (project 2014-2020.4.05.19-0001).}%
  \email{0000-0003-3052-465X}%
  \institute{Tallinn University of Technology}%
  \email{fosco.loregian@gmail.com}
  \and
  Daniele Palombi%
  \email{0000-0002-8107-5439}%
  \institute{Sapienza University of Rome}%
  \email{danielepalombi@protonmail.com}
}
\def\titlerunning{A Categorical Semantics for Hierarchical Petri Nets}
\def\authorrunning{Genovese, Herold, Loregian, Palombi}
\begin{document}
\maketitle

\begin{abstract}
  We show how a particular variety of hierarchical nets, where the firing of a transition
  in the parent net must correspond to an execution in some child net, can be modelled
  utilizing a functorial semantics from a free category -- representing the parent net --
  to the category of sets and spans between them. This semantics can be internalized via Grothendieck construction, resulting
  in the category of executions of a Petri net representing the semantics of the overall
  hierarchical net. We conclude the paper by giving an engineering-oriented overview of
  how our model of hierarchical nets can be implemented in a transaction-based smart
  contract environment.
\end{abstract}
\section{Introduction}
  \label{sec: introduction}
This paper is the fourth instalment in a series of
works~\cite{Genovese2020, Genovese2021a, Genovese2021} devoted to describing the semantics of extensions of
Petri nets using categorical tools.

Category theory has been applied to Petri nets
starting in the nineties~\cite{Meseguer1990}; see also \cite{baldan2015modular,baldan2005relating,baldan2001compositional,baldan2011mpath2pn,baldan2018petri,baldan2010petri,baldan2014encoding,baldan2008open,baldan2019petri,baldan2015asynchronous}. The main idea is that we can
use different varieties of free monoidal categories to describe
the executions (or runs) of a net~\cite{Master2020,Genovese2019c}. These works have been
influential since they opened up an avenue of applying high-level methods
to studying Petri nets and their properties. For instance,
in~\cite{Baez2020a} the categorical approach allowed to describe glueing of nets
leveraging on colimits and double categories, while category-theory
libraries such as~\cite{Genovese2019d} can be leveraged to implement nets in a formally verified way.
These libraries implement category theory \emph{directly}, so that one could translate the
categorical definitions defining some model object directly and obtain an implementation.

In~\cite{Genovese2020}, we started another line of research, where we were able
to define a categorical semantics for coloured nets employing monoidal
functors. The Grothendieck construction was then used to internalize
this semantics, obtaining the well-known result that coloured nets can be
``compiled back'' to Petri nets.

In~\cite{Genovese2021a, Genovese2021}, we extended these ideas further, and we were able to
characterize bounded nets and mana-nets -- a new kind of nets useful
to model chemical reactions -- in terms of generalized functorial semantics.

This approach, based on the correspondence between slice categories and
lax monoidal functors to the category of spans ~\cite{Pavlovic1997}, has still a lot to give. In this
paper, we show how it can be used to model hierarchical nets.

There are a lot of different ways to define hierarchical nets \cite{Jensen2009,fehling1991concept,oswald1990environment,huber1989hierarchies,Buchholz1994HierarchicalHL}, which can be seen as a graph-based model. It means that
we have one ``parent'' Petri net and a bunch of ``child'' nets. A transition firing in
the parent net corresponds to some sort of run happening in a corresponding child net.
The main net serves to orchestrate and coordinate the executions
of many child nets in the underlayer.

This paper will contain very little new mathematics. Instead, we will reinterpret
results obtained in~\cite{Genovese2020} to show how they can be used to model hierarchical nets, moreover, in a way that makes sense from an implementation perspective.

It is worth noting that category theory in this paper is used in a way that is slightly different than the usage in graph transformations research: We won't be using category theory to generalize definitions and proofs to different classes of graph(-related) objects. Instead, we will employ categorical concepts to actually build a semantics for hierarchical Petri nets.

\section{Nets and their executions}\label{sec: nets and their executions}
We start by recalling some basic constructions of category theory and some basic facts about Petri nets and their categorical formalization.
The notions of \emph{bicategory, pseudofunctor, lax functor and bimodule} are not strictly necessary to understand this paper: They only show up in results
that we cite and that could in principle be taken for granted while skimming on the details. In any case, we list these
notions here for the reader interested in parsing these results in full depth.
The first definition we recall is the one of \emph{bicategory}. Intuitively, bicategories are categories where we also allow for
``morphisms between morphisms'', called \emph{2-cells}. This in turn allows to define a version of the associativity and identity laws that is
weaker than for usual categories, holding only up to isomorphism.

\begin{definition}[Bicategory]\label{bicat}
	A \emph{(locally small) bicategory} $\mathcal{B}$ consists of the following data.
	\begin{enumerate}
		\item \label{bica:uno} A class $\mathcal{B}_o$ of \emph{objects}, denoted with Latin letters like $A,B,\dots$, also called \emph{0-cells}.
		\item \label{bica:due} A collection of (small) categories $\mathcal{B}(A,B)$, one for each $A,B\in \mathcal{B}_o$, whose objects are called \emph{1-cells} or \emph{arrows} with \emph{domain} $A$ and \emph{codomain} $B$, and whose morphisms $\alpha : f \Rightarrow g$ are called \emph{2-cells} or \emph{transformations} with domain $f$ and codomain $g$; the composition law $\circ$ in $\mathcal{B}(A,B)$ is called \emph{vertical composition} of 2-cells.
		\item A family of \emph{compositions}
		\[
			\bullet_{\mathcal{B},ABC} : \mathcal{B}(B,C)\times\mathcal{B}(A,B) \to \mathcal{B}(A,C) : (g,f)\mapsto g\bullet f
		\]
		defined for any triple of objects $A,B,C$. This is a family of functors between hom-categories, and its action on morphisms is called \emph{horizontal composition} of natural transformations, that we denote $\alpha\bullet\beta$.
		\item \label{bica:tre} For every object $A\in  \mathcal{B}_o$ there is an arrow $\id_A\in \mathcal{B}(A,A)$.
	\end{enumerate}
	To this basic structure we add\index{associator!monoidal ---}
	\begin{enumerate}
		\item a family of invertible maps $\alpha_{fgh} : (f \bullet g) \bullet h \cong f \bullet (g \bullet h)$ natural in all its arguments $f,g,h$, which taken together form the \emph{associator} isomorphisms;
		\item a family of invertible maps $\lambda_f  : \id_B \bullet f \cong f$ and $\varrho_f : f \bullet \id_A \cong f$ natural in its component $f : A \to B$, which taken together form the \emph{left unitor} and \emph{right unitor} isomorphisms.
	\end{enumerate}
	Finally, these data are subject to the following axioms.
	\index{horizontal composition}
	\index{_aaa_boxminus@$\bullet$}
	\begin{enumerate}
		\item For every quadruple of 1-cells $f,g,h,k$ we have that the diagram
		\[
			\vcenter{\xymatrix{
			((f\bullet g)\bullet h)\bullet k \ar[d]_{\alpha_{f,g,h}\bullet k}\ar[r]^{\alpha_{fg,h,k}} & (f\bullet g)\bullet (h\bullet k) \ar[r]^{\alpha_{f,g,hk}} & f\bullet (g\bullet (h\bullet k))\\
			(f\bullet (g\bullet h))\bullet k \ar[rr]_{\alpha_{f,gh,k}} && f\bullet ((g\bullet h)\bullet k)\ar[u]_{f\bullet \alpha_{g,h,k}}
			}}
		\]
		commutes.
		\item For every pair of composable 1-cells $f,g$,
		\[
			\vcenter{\xymatrix{
			(f \bullet \id_A)\bullet g\ar[dr]_{\varrho_f\bullet g}\ar[rr]^{a_{A,\id_A,g}} && f\bullet(\id_A\bullet \, g)\ar[dl]^{f\bullet \lambda_g}\\
			& f\bullet g
			}}
		\]
		commutes.
	\end{enumerate}
\end{definition}
\begin{definition}[2-category]
	A \emph{2-category} is a bicategory where the associator and unitors are the identity natural transformations. In other words, a 2-category is precisely a bicategory where horizontal composition is strictly associative, and the identities $\id_A$ work as strict identities for the horizontal composition operation.
\end{definition}
Some sources call `2-category' what we call a bicategory, and `strict 2-category' what we call a 2-category. Something similar happens for monoidal categories: a monoidal category is called \emph{strict} if its associator and left/right unitors are identity natural transformations. This is not by chance: a (strict) monoidal category $\mathcal V$ is exactly a (strict) 2-category with a single object $\ast$ (so that the category $\mathcal V$ can be identified with the category of endomorphisms of $\ast$).
\begin{example}\leavevmode
	\begin{itemize}
		\item There is a 2-category $\Cat$ where 0-cells are small categories, and the hom categories $\Cat(C,D)$ are the categories of functors and natural transformations. Composition of functors is strictly associative and unital.
		\item There is a bicategory of \emph{profunctors}, as defined in \cite{benabou2000distributors,cattani2005profunctors} and \cite[Ch. 5]{coend-calcu}. Composition of profunctors is associative up to a canonical isomorphism.
		\item Every category $\mathcal{C}$ is trivially a 2-category by taking the 2-cells to be identities. This is sometimes called the `discrete' 2-category obtained from a category $\mathcal{C}$.
		\item There is a 2-category where 0-cells are partially ordered sets $(P,\le)$, and where the category $\mathsf{Pos}(P,Q)$ is the partially ordered set of monotone functions $f : P \to Q$ and pointwise order ($f\preceq g$ iff $\forall p.fp\le gp$ in $Q$). Composition is strictly associative and unital.
	\end{itemize}
\end{example}
\begin{remark}
The fact that for every bicategory $\mathcal B$ the maps
\[
	\boxminus_{\mathcal{B},ABC} : \mathcal{B}(B,C)\times\mathcal{B}(A,B) \to \mathcal{B}(A,C) : (g,f)\mapsto g\bullet f
\]
are functors with domain a product category entails the following identity:
\begin{quote}
	Given any diagram of 2-cells like
	\[\xymatrix{
		A\ruppertwocell^{}{\alpha}
\rlowertwocell_{}{\beta}
\ar[r] & B
\ruppertwocell^{}{\gamma}
\rlowertwocell_{}{\delta}
\ar[r] & C\\
	}\]
	we have that $(\delta\bullet\beta)\circ (\gamma\bullet\alpha) = (\delta\circ\gamma) \bullet (\beta\circ\alpha)$.
	This is usually called the \emph{interchange law} in $\mathcal B$.
\end{quote}
\end{remark}

Pseudofunctors and lax functors, defined below, are some of the most widely used notions of morphism between bicategories. These are
useful to parse the deep results on which~\cref{sec: internalization} relies.

\begin{definition}[Pseudofunctor, (co)lax functor]\label{colaxe}\index{functor!pseudo---}
	Let $\mathcal{B},\mathcal{C}$ be two bicategories; a \emph{pseudofunctor} consists of
	\begin{enumerate}
		\item a function $F_o : \mathcal{B}_o \to \mathcal C_o$,
		\item a family of functors $F_{AB} : \mathcal{B}(A,B) \to \mathcal{C}(FA, FB)$,
		\item an invertible 2-cell $\mu_{fg} : Ff \circ Fg \Rightarrow F(fg)$ for each $A \xrightarrow{g}B\xrightarrow{f} C$, natural in $f$ (with respect to vertical composition) and an invertible 2-cell $\eta : \eta_f : \id_{FA} \Rightarrow F(\id_A)$, also natural in $f$.
	\end{enumerate}
	These data are subject to the following commutativity conditions for every 1-cell $A \to B$:
	\[
		\xymatrix{
		Ff\circ \id_A \ar[r]^{\varrho_{Ff}}\ar[d]_{Ff * \eta} & Ff\ar[d]^{F(\varrho_f)} & \id_B \circ Ff\ar[d]_{\eta * Ff} \ar[r]^{\lambda_{Ff}}& Ff\ar[d]^{F(\lambda_f)}\\
		Ff \circ F(\id_A)\ar[r]_{\mu_{f,\id_A}} & F(f \circ \id_A) & F(\id_B)\circ Ff\ar[r]_{\mu_{\id_B,f}} & F(\id_B \circ f)\\
		(Ff\circ Fg) \circ Fh \ar[rrr]^{\alpha_{Ff,Fg,Fh}}\ar[d]_{\mu_{fg} * Fh} &&& Ff\circ (Fg\circ Fh)\ar[d]^{Ff * \mu_{gh}}\\
		F(fg)\circ Fh \ar[d]_{\mu_{fg} * Fh} &&& Ff\circ F(gh)\ar[d]^{\mu_{f,gh}}\\
		F((fg)h) \ar[rrr]_{F \alpha_{fgh}}&&& F(f(gh))
		}
	\]
	(we denote invariably $\alpha,\lambda,\varrho$ the associator and unitor of $\mathcal{B},\mathcal{C}$).

	A \emph{lax} functor is defined by the same data, but both the 2-cells $\mu : Ff \circ Fg \Rightarrow F(fg)$ and $\eta : \id_{FA} \Rightarrow F(\id_A)$ can be non-invertible; the same coherence diagrams in Definition \ref{colaxe} hold. A \emph{colax} functor reverses the direction of the cells $\mu,\eta$, and the commutativity of the diagrams in  Definition \ref{colaxe} changes accordingly.\index{functor!lax ---}
\end{definition}

Another notion that we will make heavy use of is the one of \emph{comonad}. On the other hand, \emph{monads} and morphisms between them,
called \emph{bimodules}, will only appear in~\cref{thm: zanasi}. We will only use a straightforward consequence of this theorem, and the hurrying reader may not
linger too much on these definitions.
\begin{definition}[Monad, comonad]\label{muso_da_mona}
	Let $\CategoryC$ be a category; a \emph{monad} on $\CategoryC$ consists of an endofunctor $T : \CategoryC\to \CategoryC$ endowed with two natural transformations
	\begin{itemize}
		\item $\mu : T\circ T\Rightarrow T$, the \emph{multiplication} of the monad, and
		\item $\eta : \id_\CategoryC \Rightarrow T$, the \emph{unit} of the monad,
	\end{itemize}
	such that the following axioms are satisfied:
	\begin{itemize}
		\item the multiplication is associative, i.e. the diagram
					\[
						\vcenter{\xymatrix{
								T\circ T\circ T\ar[r]^{T *\mu}\ar[d]_{\mu *T} & T\circ T \ar[d]^\mu\\
								T\circ T \ar[r]_\mu & T
							}	}
					\]
					is commutative, i.e. the equality of natural transformations $\mu\circ (\mu * T) = \mu \circ (T * \mu)$ holds;
		\item the multiplication has the transformation $\eta$ as unit, i.e. the diagram
					\[
						\vcenter{\xymatrix{
								T \ar[r]^{\eta *T}\ar@{=}[dr]& T\circ T \ar[d]_\mu & T\ar[l]_{T*\eta} \ar@{=}[dl]\\
								& T &
							}	}
					\]
					is commutative, i.e. the equality of natural transformations $\mu\circ (\eta *T)=\mu\circ (T * \eta)= \id_T$ holds.
	\end{itemize}
	Dually, let $\CategoryC$ be a category; a \emph{comonad} on $\CategoryC$ consists of an endofunctor $T : \CategoryC\to \CategoryC$ endowed with two natural transformations
	\begin{itemize}
		\item $\sigma : T \Rightarrow T\circ T$, the \emph{comultiplication} of the comonad, and
		\item $\epsilon : T\Rightarrow \id_\CategoryC$, the \emph{counit} of the comonad,
	\end{itemize}
	such that the following axioms are satisfied:
	\begin{itemize}
		\item the comultiplication is coassociative, i.e. the diagram
					\[
						\vcenter{\xymatrix{
								T\circ T\circ T\ar@{<-}[r]^{T *\sigma}\ar@{<-}[d]_{\sigma *T} & T\circ T \ar@{<-}[d]^\sigma\\
								T\circ T \ar@{<-}[r]_\sigma & T
							}	}
					\]
					is commutative.
		\item the comultiplication has the transformation $\epsilon$ as counit, i.e. the diagram
					\[
						\vcenter{\xymatrix{
								T \ar@{<-}[r]^{\epsilon *T}\ar@{=}[dr]& T\circ T \ar@{<-}[d]_\sigma & T\ar@{<-}[l]_{T*\epsilon} \ar@{=}[dl]\\
								& T &
							}	}
					\]
					is commutative.
	\end{itemize}
\end{definition}
\begin{definition}[Bimodule]
	Given a bicategory $\mathcal B$ having finite colimits (in the 2-categorical sense of \cite{2catlimits}), define the 2-category $\cate{Mod}(\mathcal B)$ of \emph{bimodules} as in \cite[2.19]{Zanasi2018}:
	\begin{itemize}
		\item 0-cells are the monads in $\mathcal B$;
		\item 1-cells $T \to S$ are \emph{bimodules}, i.e. 1-cells $H : C \to D$ (assuming $T$ is a monad on $C$, and $S$ a monad on $D$) equipped with suitable action maps: $\rho : HT \to H$ and $\lambda : SH\to H$ satisfying suitable axioms expressing the fact that $T$ acts on the right over $H$, via $\rho$ (resp., $S$ acts on the left on $H$, va $\lambda$);
		\item 2-cells are natural transformations $\alpha : H \Rightarrow K :  T\to S$ compatible with the action maps.
	\end{itemize}
	\end{definition}
\subsection{Categorical Petri nets}
Having recalled some of the category theory we are going to use, we now summarize some needed definitions underlying the
study of Petri nets from a categorical perspective.

\begin{figure}
  \begin{small}
    \begin{center}
    \begin{adju}
	\begin{tikzpicture}
		\pgfmathsetmacro\bS{5}
		\pgfmathsetmacro\hkX{(\bS/3.5)}
		\pgfmathsetmacro\kY{-1.5}
		\pgfmathsetmacro\hkY{\kY*0.5}

		\draw pic (m0) at (0,0) {netA={{1}/{1}/{2}/{0}/{}/{}/{}}};
		\draw pic (m1) at (\bS,0) {netA={{0}/{2}/{2}/{0}/{$\blacktriangle$}/{}/{}}};
		\draw pic (m2) at ({2 * \bS},0) {netA={{0}/{1}/{3}/{1}/{}/{$\blacktriangledown$}/{}}};
		\draw pic (m3) at ({3 * \bS},0) {netA={{0}/{1}/{2}/{2}/{}/{}/{$\blacktriangledown$}}};

		\begin{scope}[very thin]
		\foreach \j in {1,...,3} {
			\pgfmathsetmacro \k { \j * \bS - 1 };
			\draw[gray,dashed] (\k,-4) -- (\k,-8.25);
			\draw[gray] (\k,1) -- (\k,-4);
		}
		\end{scope}

		\begin{scope}[shift={(0,-4)}, oriented WD, bbx = 1cm, bby =.4cm, bb min width=1cm, bb port sep=1.5]

		\draw node [fill=\backgrnd,bb={1}{1}] (Tau) at (\bS -1,-1) {$t$};
		\draw node [fill=\backgrnd,bb={1}{2}, ] (Mu)  at ({2 * \bS - 1},-1) {$v$};
		\draw node [fill=\backgrnd,bb={1}{1}] (Nu)  at ({3 * \bS - 1},{2 * \kY}) {$u$};

		\draw (-1,-1) --     node[above] {$p_1$}       (0,-1)
		--                  node[above] {}          (Tau_in1);

		\draw (-1,-2) -- node[above] {$p_2$} (0,-2) -- (\bS-1, -2);
		\draw (-1,-3) -- node[above] {$p_3$} (0,-3) -- (\bS-1, -3);
		\draw (-1,-4) -- node[above] {$p_3$} (0,-4) -- (\bS-1, -4);

		\draw (Tau_out1) -- node[above] {$p_2$}    (Mu_in1);
		\draw (\bS-1,-2) -- (2*\bS-1, -2);
		\draw (\bS-1,-3) -- (2*\bS-1, -3);
		\draw (\bS-1,-4) -- (2*\bS-1, -4);

		\draw (Mu_out1) -- node[above] {$p_3$}    (3*\bS-1, -0.725);
		\draw (Mu_out2) -- node[above] {$p_4$}    (3*\bS-1, -1.325);

		\draw (2*\bS-1,-2) -- (3*\bS-1, -2);
		\draw (2*\bS-1,-3) -- (Nu_in1);
		\draw (2*\bS-1,-4) -- (3*\bS-1, -4);

		\draw (3*\bS-1,-0.725) to (4*\bS-2, -1.325) -- node[above] {$p_3$} (4*\bS-1, -1.325);
		\draw (3*\bS-1,-1.325) -- (3*\bS,-1.325) to (4*\bS-2, -4) -- node[above] {$p_4$} (4*\bS-1, -4);
		\draw (3*\bS-1,-2) to (4*\bS-2, -0.725) -- node[above] {$p_2$} (4*\bS-1, -0.725);
		\draw (Nu_out1) to (4*\bS-2, -3) -- node[above] {$p_4$} (4*\bS-1, -3);
		\draw (3*\bS-1,-4) to (4*\bS-2, -2) -- node[above] {$p_3$} (4*\bS-1, -2);

		\end{scope}

		\begin{pgfonlayer}{background}
		\filldraw [line width=4mm,join=round,\backgrnd]
		((-1,1.5)  rectangle (19,-8.5);
		\end{pgfonlayer}

		\end{tikzpicture}
	\end{adju}
    \end{center}
  \caption{Graphical representation of a net's execution.}\label{fig: execution of a net}
  \end{small}
\end{figure}
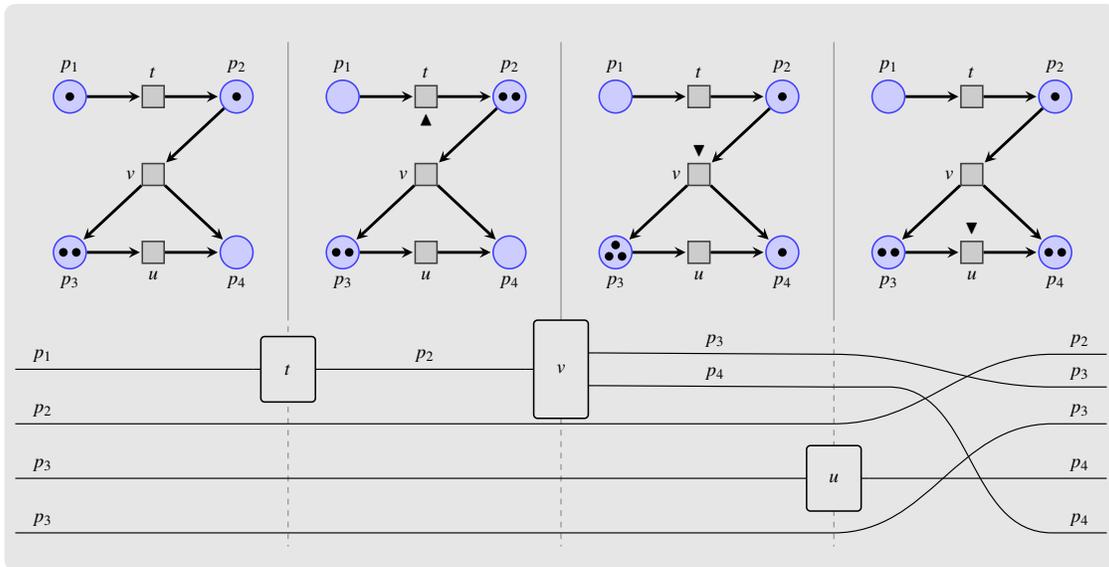

\begin{notation}
  Let $S$ be a set; a multiset is a function $S \to \Naturals$.
  Denote with $\Msets{S}$ the set of multisets over $S$.
  Multiset sum and difference (only partially defined)
  are defined pointwise and will be denoted with $\oplus$ and
  $\ominus$, respectively.
  The set $\Msets{S}$ together with $\oplus$ and the empty multiset is isomorphic to the free commutative monoid on $S$.
\end{notation}
\begin{definition}[Petri net]\label{def: Petri net}
  A
  \emph{Petri net} is a pair of functions $T \xrightarrow{s,t} \Msets{S}$ for
  some sets $T$ and $S$, called the set of places and transitions of the net, respectively.
  $s,t$ are called \emph{input and output} functions, respectively, or equivalently \emph{source
  and target}.

  A \emph{morphism of nets} is a pair of functions $f: T \to T'$ and $g: S \to S'$ such that
  the following square commutes, with $\Msets{g}: \Msets{S} \to \Msets{S'}$ the obvious
  lifting of $g$ to multisets:
  \begin{equation*}
    \begin{tikzcd}
      {\Msets{S}} & {T} & {\Msets{S}} \\
      {\Msets{S'}} & {T'} & {\Msets{S'}}
      \arrow["{s}"', from=1-2, to=1-1]
      \arrow["{s'}", from=2-2, to=2-1]
      \arrow["{t'}"', from=2-2, to=2-3]
      \arrow["{t}", from=1-2, to=1-3]
      \arrow["{\Msets{g}}"', from=1-1, to=2-1]
      \arrow["{\Msets{g}}", from=1-3, to=2-3]
      \arrow["{f}" description, from=1-2, to=2-2]
    \end{tikzcd}
  \end{equation*}
  Petri nets and their morphisms form a category, denoted $\Petri$. Details can be found in~\cite{Meseguer1990}.
\end{definition}
\begin{definition}[Markings and firings]\label{def: Petri net firing}
  A \emph{marking} for a net  $T \xrightarrow{s,t} \Msets{S}$ is an element of $\Msets{S}$,
  representing a distribution of tokens in the net places. A transition $u$ is \emph{enabled} in a marking $M$ if
  $M \ominus s(u)$ is defined. An enabled transition can \emph{fire}, moving tokens in the net.
  Firing is considered an atomic event, and the marking resulting from firing $u$ in $M$ is $M \ominus s(u) \oplus t(u)$.
  Sequences of firings are called \emph{executions}.
\end{definition}
The main insight of categorical semantics for Petri nets is that the information contained in
a given net is enough to generate a free symmetric strict monoidal category representing all the
possible ways to run the net. There are multiple ways to do
this~\cite{Sassone1995,Genovese2019c,Genovese2019b,Master2020, Baez2021}. In this work,
we embrace the \emph{individual-token philosophy},
where tokens are considered distinct and distinguishable and thus require
the category in~\cref{def: executions individual token philosophy} to
have non-trivial symmetries.
\begin{definition}[Category of executions -- individual-token philosophy]\label{def: executions individual token philosophy}
  Let $N: T \xrightarrow{s,t} \Msets{S}$ be a Petri net.
  We can generate a \emph{free symmetric strict monoidal category (\fssmc)}, $\Free{N}$, as follows:
  \begin{itemize}
    \item The monoid of objects is the free monoid generated by $S$. Monoidal product of objects $A,B$ is denoted with $A \otimes B$.
    \item Morphisms are generated by $T$: each $u \in T$ corresponds to a morphism generator
    $(u,su,tu)$, pictorially represented as an arrow $su \xrightarrow{u} tu$;
    morphisms are obtained by considering all the formal (monoidal) compositions of generators and identities.
  \end{itemize}
  A detailed description of this construction can be found in~\cite{Master2020}.
\end{definition}
In this definition, objects represent markings of a net. For instance, the object $A\oplus A \oplus B$ means ``two tokens
in $A$ and one token in $B$''. Morphisms represent executions of a net, mapping markings to markings.
A marking is reachable from another one if and only if there is a morphism between them. An example is provided in~\cref{fig: execution of a net}.

\section{Hierarchical nets}\label{sec: hierarchical nets}
Now we introduce the main object of study of the paper,
\emph{hierarchical nets}. As we pointed out in~\cref{sec: introduction},
there are many different ways to model hierarchy in Petri nets~\cite{Jensen2009},
often incompatible with each other.
We approach the problem from a developer's perspective,
wanting to model the idea that ``firing a transition'' amounts to call
another process and waiting for it to finish. This is akin to calling
subroutines in a piece of code. Moreover, we do not want to destroy the
decidability of the reachability relation for our nets~\cite{Esparza1994}, as it happens for
other hierarchical models such as the net-within-nets
framework~\cite{Kohler-Bussmeier2014}.
We consider this to be an essential requirement
for practical reasons.

We will postpone any formal definition to~\cref{sec: semantics for hierarchical nets}. In the present work, we focus on
giving an intuitive explanation of what our requirements are.
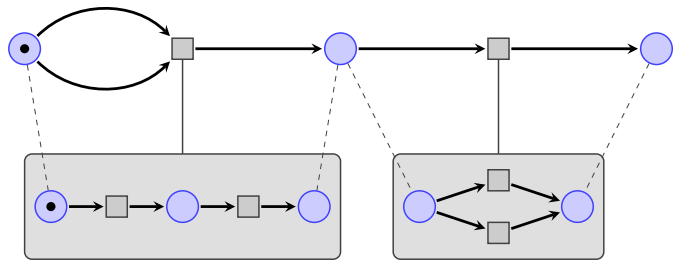
\begin{figure}[!ht]\centering
  \scalebox{0.7}{
  \begin{tikzpicture}[node distance=1.3cm,>=stealth',bend angle=45,auto]
    %
    \filldraw[thick, draw=black!75, fill=gray!25, rounded corners] (0,-1) rectangle (6,-3);
    \draw[black!75, thick]  (3,1) -- (3,-1);
    \filldraw[thick, draw=black!75, fill=gray!25, rounded corners] (7,-1) rectangle (11,-3);
    \draw[black!75, thick]  (9,1) -- (9,-1);
    \node [place,tokens=1] (1a) at (0,1){};
    \node [transition] (2a) at (3,1)     {}
        edge [inarrow, pre, bend left] (1a)
        edge [inarrow, pre, bend right] (1a);
    \node [place,tokens=0] (3a) at (6,1)  {}
        edge [inarrow, pre] (2a);
    \node [transition] (4a) at (9,1) {}
        edge [inarrow, pre] (3a);
    \node [place,tokens=0] (5a) at (12,1) {}
        edge [inarrow, pre] (4a);
    %
    %
    %
    \node [place,tokens=1] (c11a) at (0.5,-2){};
    \node [transition] (c12a) at (1.75,-2)     {}
        edge [inarrow, pre] (c11a);
    \node [place,tokens=0] (c13a) at (3,-2)  {}
        edge [inarrow, pre] (c12a);
    \node [transition] (c14a) at (4.25,-2)     {}
        edge [inarrow, pre] (c13a);
    \node [place,tokens=0] (c15a) at (5.5,-2)  {}
        edge [inarrow, pre] (c14a);
    \draw[black!74, dashed] (1a) -- (c11a);
    \draw[black!74, dashed] (3a) -- (c15a);
    %
    %
    %
    \node [place,tokens=0] (c23a) at (7.5,-2)  {};
    \node [transition] (c24a) at (9,-1.5) {}
        edge [inarrow, pre] (c23a);
    \node [transition] (c24b) at (9,-2.5) {}
    edge [inarrow, pre] (c23a);
    \node [place,tokens=0] (c25a) at (10.5,-2) {}
        edge [inarrow, pre] (c24a)
        edge [inarrow, pre] (c24b);
    \draw[black!74, dashed] (3a) -- (c23a);
    \draw[black!74, dashed] (5a) -- (c25a);
  \end{tikzpicture}}
  \caption{A hierarchical net.}\label{fig: hierarchical net}
\end{figure}

Looking at the net in~\cref{fig: hierarchical net}, we see a
net on the top, which we call \emph{parent}. To each transition
of the parent net is attached another net, which we call \emph{child}. 
Transitions can only have one child, but the parent net may have 
multiple transitions, and hence multiple children overall.
Connecting input and output places of a transition in the parent
net with certain places in the corresponding child, we can represent
the orchestration by saying that each time a transition in the parent net
fires, its input tokens are transferred to the corresponding child net,
that takes them across until they reach a place connected with the
output place in the parent net. This way, the atomic act of firing a transition
in the parent net results in an execution of the corresponding child.
\begin{figure}[!ht]\centering
  \scalebox{0.7}{
  \begin{tikzpicture}[node distance=1.3cm,>=stealth',bend angle=45,auto]
    %
    \node [place,tokens=1] (c11a) at (0.5,-2){};
    \node [transition] (c12a) at (1.75,-2)     {}
        edge [inarrow, pre] (c11a);
    \node [place,tokens=0] (c13a) at (3,-2)  {}
        edge [inarrow, pre] (c12a);
    \node [transition] (c14a) at (4.25,-2)     {}
        edge [inarrow, pre] (c13a);
    \node [place,tokens=0] (c15a) at (6,-2)  {}
        edge [inarrow, pre] (c14a);
    %
    %
    %
    \node [place,tokens=0] (c23a) at (6,-2)  {};
    \node [transition] (c24a) at (9,-1.5) {}
        edge [inarrow, pre] (c23a);
    \node [transition] (c24b) at (9,-2.5) {}
    edge [inarrow, pre] (c23a);
    \node [place,tokens=0] (c25a) at (10.5,-2) {}
        edge [inarrow, pre] (c24a)
        edge [inarrow, pre] (c24b);
  \end{tikzpicture}}
  \caption{Replacing transitions in the parent net of~\cref{fig: hierarchical net} with its children.}\label{fig: replacing transitions in the parent net with child nets}
\end{figure}
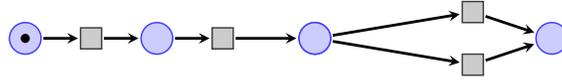

Notice that we are not interested in considering the semantics
of such hierarchical net to be akin to the one in~\cref{fig: replacing transitions in the parent net with child nets},
where we replaced transitions in the parent net with their corresponding
children. Indeed, this way of doing things is similar to what happens in~\cite{oswald1990environment}: 
In this model, transitions in the parent net are considered as placeholders for the children nets.
There are two reasons why we distance ourselves from this approach: First, we want to consider
transition firings in the parent net as atomic events, and replacing
nets as above destroys this property. Secondly, such replacement
is not so conceptually easy given that we do not
impose any relationship between the parent
net's topologies and its children. Indeed, the leftmost transition of the parent net
in~\cref{fig: hierarchical net} consumes two inputs, while the corresponding
leftmost transition in its child only takes one. How do we account for this
in specifying rewriting-based semantics for hierarchical nets?

\section{Local semantics for Petri nets}\label{sec: local semantics for Petri nets}
We concluded the last section pointing out reasons that
make defining a semantics for hierarchical nets less intuitive
than one would initially expect. Moreover, requiring the transition firings
in the parent net to be considered as atomic events basically rules out
the majority of the previous approaches to hierarchical Petri nets, as the one sketched 
in~\cite{Jensen2009, oswald1990environment}.
Embracing an engineering perspective, we could get away with some
ad-hoc solution to conciliate that parent and child net topologies
are unrelated. One possible way, for instance, would be imposing constraints
between the shapes of the parent net and its children. However, in
defining things ad-hoc, the possibility for
unforeseen corner cases and situations we do not know how to deal with
becomes high. To avoid this, we embrace a categorical perspective
and define things up to some degree of canonicity.

Making good use of the categorical work already carried out on
Petri nets, our goal is to leverage it and get to a plausible definition of categorical
semantics for hierarchical nets. Our strategy is
to consider a hierarchical net as an extension of a Petri net: The
parent net will be the Petri net we extend, whereas the children
nets will be encoded in the extension.

This is precisely the main idea contained in~\cite{Genovese2020}, that is, the idea of describing net
extensions with different varieties of monoidal functors. Indeed,
we intend to show how the theory presented in~\cite{Genovese2020}, and
initially serving a wholly different purpose, can be reworked
to represent hierarchical nets with minimal effort.

As for semantics, we will use strict monoidal functors and name it
\emph{local} because the strict-monoidality requirement
amounts to endow tokens with properties
that cannot be shared with other tokens. To understand this
choice of naming a little bit better, it may be worth comparing it with
the notion of \emph{ non-local semantics}, defined in
terms of lax-monoidal-lax functors, that we gave in~\cite{Genovese2021a}.
\begin{definition}[Local semantics for Petri nets]\label{def: PetriS}
  Given a strict monoidal category $\Semantics$, a
  \emph{Petri net with a local $\Semantics$-semantics}
   is a pair  $\NetSem{N}$, consisting of a Petri net $N$ and
   a strict monoidal functor
  \[\Fun{N}: \Free{N} \to \Semantics.\]
  A morphism $F: \NetSem{M} \to \NetSem{N}$
  is just a strict monoidal functor $F: \Free{M} \to \Free{N}$ such that
  $\Fun{M} = F \cp \Fun{N}$, where we denote composition in diagrammatic order;
  i.e.\ given $f\colon c\to d$ and $g\colon d\to e$,
  we denote their composite by $(f\cp g)\colon c\to e$.

  Nets equipped with $\Semantics$-semantics and their morphisms form a monoidal category denoted $\PetriS{\Semantics}$, with the monoidal structure arising from the product in $\Cat$.
\end{definition}
In~\cite{Genovese2020}, we used local semantics to describe guarded Petri nets,
using $\Span$ as our category of choice. We briefly summarize this,
as it will become useful later.
\begin{definition}[The category $\Span$]
  We denote by $\Span$  the 1-category of sets and spans,
  where isomorphic spans are identified. This category is symmetric monoidal.
  From now on, we will work with the \emph{strictified} version of $\Span$, respectively.
\end{definition}
\begin{notation}
  Recall that a morphism $A\to B$ in $\Span$ consists of a set $S$ and a pair of functions $A\from S\to B$. When we need to extract this data from $f$, we write
  \begin{equation*}
    A\From{f_1}S_f\To{f_2}B
  \end{equation*}
  We sometimes consider the span as a function $f\colon S_f\to A\times B$,
  thus we may write $f(s)=(a,b)$ for $s\in S_f$ with $f_1(s)=a$ and $f_2(s)=b$.
\end{notation}
\begin{definition}[Guarded nets with side effects]
  \label{def: guarded net}
  A \emph{guarded net with side effects} is an object
  of $\PetriSpan$. A morphism of guarded nets with side effects
  is a morphism in $\PetriSpan$.
\end{definition}
\begin{example}
  Let us provide some intuition behind the definition of $\PetriSpan$.

  Given a net $N$, its places (generating objects of $\Free{N}$) are sent to sets. Transitions (generating morphisms of $\Free{N}$) are mapped to spans.
  Spans can be understood as \emph{relations with witnesses}, provided by elements in the apex of the span: Each path from the span domain to its codomain is indexed by some element of the span apex, as it is shown in \cref{fig: semantics in span}. Witnesses allow considering different paths between the same elements. These paths represent the actions of processing the property a token is enowed with according to some side effect. Indeed,
  an element in the domain can be sent to different elements in the codomain via different paths. We interpret this as \emph{non-determinism}:
  the firing of the transition is not only a matter of the tokens input and output; it also includes the chosen path, which we interpret as having side-effects interpreted outside of our model.
\end{example}
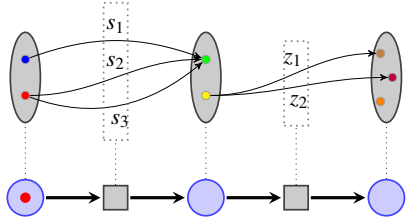
\begin{figure}[!ht]\centering
    \scalebox{0.8}{
    \begin{tikzpicture}[node distance=1.3cm,>=stealth',bend angle=45,auto]
      \node [place,colored tokens={red},] (1a) at (0,0){};
      \node [transition] (2a) at (1.5,0)     {}
          edge [inarrow, pre] (1a);
      \node [place,tokens=0] (3a) at (3,0)  {}
          edge [inarrow, pre] (2a);
      \node [transition] (4a) at (4.5,0) {}
          edge [inarrow, pre] (3a);
      \node [place,tokens=0] (5a) at (6,0) {}
          edge [inarrow, pre] (4a);
      \draw[thick, draw=black!75, fill=black!20] (0,2) ellipse (0.25cm and 0.75cm);
      \draw[thick, draw=black!75, fill=black!20] (3,2) ellipse (0.25cm and 0.75cm);
      \draw[thick, draw=black!75, fill=black!20] (6,2) ellipse (0.25cm and 0.75cm);

      \node[circle,fill=red, draw=gray!75, inner sep=0pt,minimum size=4pt] (red) at (0,1.7) {};
      \node[circle,fill=blue, draw=gray!75, inner sep=0pt,minimum size=4pt] (blue) at (0,2.3) {};

      \node[circle,fill=yellow, draw=gray!75, inner sep=0pt,minimum size=4pt] (yellow) at (3,1.7) {};
      \node[circle,fill=green, draw=gray!75, inner sep=0pt,minimum size=4pt] (green) at (3,2.3) {};

      \node[circle,fill=brown, draw=gray!75, inner sep=0pt,minimum size=4pt] (brown) at (5.9,2.4) {};
      \node[circle,fill=purple, draw=gray!75, inner sep=0pt,minimum size=4pt] (purple) at (6.1,2) {};
      \node[circle,fill=orange, draw=gray!75, inner sep=0pt,minimum size=4pt] (orange) at (5.9,1.6) {};

      \draw[->, out=20, in=160] (blue) to node[midway, above] {$s_1$} (green);
      \draw[->, out=0, in=180] (red) to node[midway, above] {$s_2$} (green);
      \draw[->, out=-20, in=220] (red) to node[midway, below] {$s_3$} (green);
      \draw[->, out=0, in=180] (yellow) to node[midway, below] {$z_2$} (purple);
      \draw[->, out=0, in=180] (yellow) to node[midway, above] {$z_1$} (brown);

      \draw[dotted, -]  (1a.north) -- (0,1.25);
      \draw[dotted, -]  (3a.north) -- (3,1.25);
      \draw[dotted, -]  (5a.north) -- (6,1.25);

      \node[transition,dotted,fill=none, draw=black!50, inner sep=0pt,minimum height=65pt] (left) at (1.5,2.1) {};
      \node[transition,dotted,fill=none, draw=black!50, inner sep=0pt,minimum height=40pt] (right) at (4.5,1.9) {};

      \draw[dotted, -]  (2a.north) -- (left.south);
      \draw[dotted, -]  (4a.north) -- (right.south);
    \end{tikzpicture}}
    \caption{Semantics in $\Span$}\label{fig: semantics in span}
\end{figure}
In \cref{fig: semantics in span} the composition
of paths is the empty span: Seeing things from
a reachability point of view, the process given by firing the left transition
and then the right will never occur. This is because the rightmost transition has a guard that only accepts yellow tokens, so that a green token can never be processed by it. This is witnessed by the fact that there is no path connecting the green dot with any dot on its right. The relation with reachability can be made precise by
recasting~\cref{def: executions individual token philosophy}.
\begin{definition}[Markings for guarded nets]\label{def: marking guarded Petri net}
  Given a guarded Petri net with side effects $\NetSem{N}$, a \emph{marking}
  for $\NetSem{N}$ is a pair $(X, x)$ where $X$ is an object of $\Free{N}$ and
  $x \in \Fun{N}X$. We say that a marking $(Y,y)$ is \emph{reachable} from $(X,x)$
  if there is a morphism $f: X \to Y$ in $\Free{N}$ and an element $s \in S_f$ such that $\Fun{N}f(s) = (x,y)$.
\end{definition}

\section{Semantics for hierarchical nets}\label{sec: semantics for hierarchical nets}
In the span semantics we
can encode externalities in the tips of the spans to which we send transitions.
That is, given a bunch of tokens endowed with some properties,
to fire a transition, we need to provide a \emph{witness} that testifies how
these properties have to be handled.
The central intuition of this paper is that we can use side effects to encode
the runs of some other net: To fire a transition in the parent net, we need to provide
a \emph{trace} of the corresponding child net. So we are saying that to fire a transition in the
parent net, a valid execution of the corresponding child net must be provided. Relying on the results in \cref{sec: nets and their executions},
we know that such valid executions are exactly the morphims in the free symmetric strict monoidal category generated by the child net.
Putting everything together, we want the tips of our spans to ``represent'' morphisms in the monoidal
categories corresponding to the children nets.
The following result makes this intuition precise, explaining how monoidal categories and spans are related:
\begin{theorem}[{\cite[Section 2.4.3]{Zanasi2018}}]\label{thm: zanasi}
  Given a category $A$ with finite limits, a \emph{category internal in $A$} is a monad in $\Span(A)$.
  Categories are monads in $\Span$, whereas strict monoidal categories are monads in $\Span(\Mon)$, with
  $\Mon$ being the category of monoids and monoid homomorphisms. A symmetric monoidal category is
  a \emph{bimodule} in $\Span(\Mon)$.
\end{theorem}
It is worth pointing out, at least intuitively, how this result works: Given a category $\CategoryC$, we denote with
$\Object{\CategoryC}$ and $\Arrow{\CategoryC}$ the sets\footnote{Here we are assuming that the objects and morphisms of our categories aren't proper classes.
This assumption is harmless in our context unless one wants to consider a Petri net whose places and transitions, respectively, form a proper class.}
of objects and morphisms of $\CategoryC$, respectively. Then we can form a span:
\begin{equation*}
  \Object{\CategoryC} \xleftarrow{\Dom{}} \Arrow{\CategoryC} \xrightarrow{\Cod{}} \Object{\CategoryC}
\end{equation*}
where the legs send a morphism to its domain and codomain, respectively. This is clearly not enough,
since in a category we have a notion of identity and composition, but asking for a monad provides exactly this.
For instance, the monad multiplication in this setting becomes a span morphism
\begin{equation*}
  \scalebox{0.89}{
   \begin{tikzpicture}
    \node (tip) at (0,3) {$\Arrow{\CategoryC} \times_{\Object{\CategoryC}} \Arrow{\CategoryC}$};
    \node (tipl) at (-2,1.5) {$\Arrow{\CategoryC}$};
    \node (tipr) at (2,1.5) {$\Arrow{\CategoryC}$};
    \node (int) at (-4,0) {$\Object{\CategoryC}$};
    \node (net) at (0,0) {$\Object{\CategoryC}$};
    \node (outt) at (4,0) {$\Object{\CategoryC}$};

    \begin{scope}[>={Stealth[black]},
      every node/.style={fill=white,circle},
      every edge/.style={draw=black}]

      \draw[->] (tip) to (tipl);
      \draw[->] (tip) to (tipr);
      \draw[->] (tipl) to node[] {$\Dom{}$} (int);
      \draw[->] (tipl) to node[] {$\Cod{}$} (net);
      \draw[->] (tipr) to node[] {$\Dom{}$} (net);
      \draw[->] (tipr) to node[] {$\Cod{}$} (outt);

      \path (tip) to node[near start, rotate=45] {$\urcorner$} (net);

      \node(explanation) at (6,1.5) {$\xrightarrow{\text{monad multiplication}}$};
      \node (cat) at (10,1.5) {$\Arrow{\CategoryC}$};
      \node (catl) at (8,1.5) {$\Object{\CategoryC}$};
      \node (catr) at (12,1.5) {$\Object{\CategoryC}$};
      \draw[->] (cat) to node[above] {$\Dom{}$} (catl);
      \draw[->] (cat) to node[above] {$\Cod{}$} (catr);
    \end{scope}
  \end{tikzpicture}}
  \end{equation*}
which gives composition of arrows. Similarly, the monad unit singles out identities,
and the monad laws witness the associativity and identity laws. In a similar way, monoidal categories are represented as above, but we
furthermore require $\Object{\CategoryC}$ and $\Arrow{\CategoryC}$ to be endowed with a monoid structure (representing the action of the monoidal
structure on the objects and morphisms of $\CategoryC$, respectively), and that this structure is preserved by the span legs, while the bimodule structure on
top of the monad witnesses the monoidal symmetries.

For the scope of our applications, we remember that each Petri net $N$ generates a free symmetric strict monoidal category $\Free{N}$,
which will correspond to a bimodule in $\Span(\Mon)$. So, in particular, we have a span of monoids\footnote{We are abusing notation,
and writing $\Object{N}$, $\Arrow{N}$ in place of $\Object{\Free{N}}$, $\Arrow{\Free{N}}$, respectively.}
\begin{equation*}
\Object{N} \xleftarrow{\Dom{}} \Arrow{N} \xrightarrow{\Cod{}} \Object{N}
\end{equation*}
underlying a bimodule, with $\Object{N}$ and $\Arrow{N}$, representing the objects
and arrows of the category, respectively, both free. We will refer to such a span as
\emph{the \fssmc $N$ (in $\Span(\Mon)$)}.
\begin{definition}[Hierarchical nets -- External definition]\label{def: hierarchical nets external definition}
  A \emph{hierarchical net} is a functor $\Free{N} \to \Span(\Mon)$ defined as follows:
  \begin{itemize}
    \item Each generating object $A$ of $\Free{N}$ is sent to a set $FA$, aka
    \emph{the set of accepting states for the place $A$}.
    \item Each generating morphism $A \xrightarrow{f} B$ is sent to a span with the following
    shape:
    \begin{equation*}
    \scalebox{0.89}
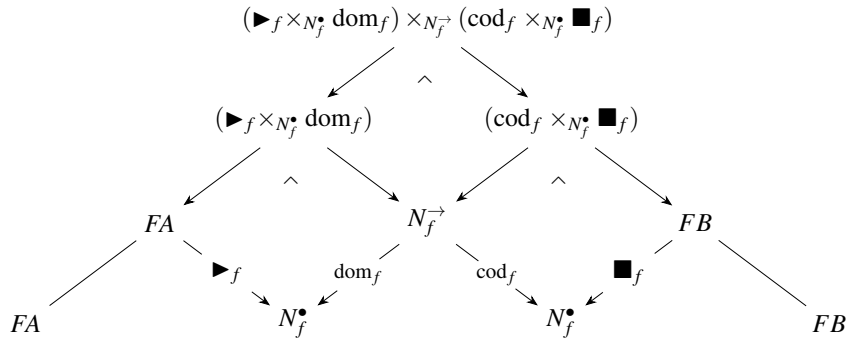
{
     \begin{tikzpicture}
      \node (tip) at (0,4.5) {$(\Play{f} \mathrel{\times_{\Object{N_f}}} \Dom{f}) \times_{\Arrow{N_f}} (\Cod{f} \mathrel{\times_{\Object{N_f}}} \Stop{f})$};
      \node (tipl) at (-2,3) {$(\Play{f} \mathrel{\times_{\Object{N_f}}} \Dom{f})$};
      \node (tipr) at (2,3) {$(\Cod{f} \mathrel{\times_{\Object{N_f}}} \Stop{f})$};
      \node (int) at (-4,1.5) {$FA$};
      \node (net) at (0,1.5) {$\Arrow{N_f}$};
      \node (outt) at (4,1.5) {$FB$};
      \node (inl) at (-6,0) {$FA$};
      \node (inr) at (-2,0) {$\Object{N_f}$};
      \node (outl) at (2,0) {$\Object{N_f}$};
      \node (outr) at (6,0) {$FB$};

      \begin{scope}[>={Stealth[black]},
        every node/.style={fill=white,circle},
        every edge/.style={draw=black}]

        \draw[->] (tip) to (tipl);
        \draw[->] (tip) to (tipr);
        \draw[->] (tipl) to (int);
        \draw[->] (tipl) to (net);
        \draw[->] (tipr) to (net);
        \draw[->] (tipr) to (outt);

        \draw[-] (int) to (inl);
        \draw[->] (int) to node[midway] {$\Play{f}$} (inr);
        \draw[->] (net) to node[midway, inner sep=0pt] {\footnotesize $\Dom{f}$} (inr);
        \draw[->] (net) to node[midway, inner sep=0pt] {\footnotesize $\Cod{f}$} (outl);
        \draw[->] (outt) to node[midway] {$\Stop{f}$} (outl);
        \draw[-] (outt) to (outr);

        \path (tip) to node[near start, rotate=45] {$\urcorner$} (net);
        \path (tipl) to node[near start, rotate=45] {$\urcorner$} (inr);
        \path (tipr) to node[near start, rotate=45] {$\urcorner$} (outl);

      \end{scope}
    \end{tikzpicture}}
    \end{equation*}
    The \fssmc $N_f$ at the center of the span is called the \emph{child net associated to $f$};
    the morphisms $\Play{f}$ and $\Stop{f}$ are called \emph{play $N_f$} and
    \emph{stop $N_f$}, respectively.
  \end{itemize}
\end{definition}
Unrolling the definition, we are associating to each generating morphism
of $f$ of $\Free{N}$ -- the parent net -- a \fssmc $N_f$-- the child net. As the feet of
the spans corresponding to the child nets will, in general, be varying with the net themselves, we need to pre and post-compose them with other spans
to ensure composability: $\Play{f}$ and $\Stop{f}$
represent morphisms that select the \emph{initial and accepting states} of $N_f$, that is,
markings of $N_f$ in which the computation starts, and markings of $N_f$ in which the computation
is considered as concluded.
Notice how this also solves the problems highlighted in~\cref{sec: hierarchical nets},
as $\Play{f}$ and $\Stop{f}$ mediate between the shape of inputs/outputs of the
transition $f$ and the shape of $N_f$ itself.
\begin{remark}\label{rem: interpreting hierarchical nets}
  Interpreting markings as in~\cref{def: marking guarded Petri net}, We see that to
  fire $f$ in the parent net we need to provide a triple $(a,x,b)$, where:
  \begin{itemize}
    \item $a$ is an element of $FA$, witnessing that the tokens
    in the domain of $f$ are a valid initial state for $N_f$.
    \item $x$ is an element of $\Arrow{N_f}$, that is,
    a morphism of $N_f$, and hence an execution of the child net.
    This execution starts from the marking $\Play{f}a$ and ends in the marking $\Stop{f}b$.
    \item $b$ is an element of $FB$, witnessing that the resulting
    state of the execution $x$ is \emph{accepting}, and can be lifted
    back to tokens in the codomain of $f$.
  \end{itemize}
\end{remark}
\begin{definition}[Category of hierarchical Petri nets]
  Nets $\NetSem{N}$ in the category $\PetriSpan$ with $\Fun{N}$ having the shape
  of~\cref{def: hierarchical nets external definition} form a
  subcategory, denoted with $\PetriHier$, and called \emph{the category
  of hierarchical Petri nets}.
\end{definition}
\begin{remark}
  Using the obvious forgetful functor $\Mon \to \Set$ we obtain
  a functor $\Span(\Mon) \to \Span$, which allows to recast our
  non-local semantics in a more liberal setting. In particular, we could send a transition to spans whose components are \emph{subsets} of the monoids heretofore considered. We could select only a subset of the executions/states of the child net as valid witnesses to fire a transition in the parent.

  Everything we do in this work will go through smoothly, but we consider this approach less elegant; thus, we will not mention it anymore.
\end{remark}

\section{Internalization}\label{sec: internalization}
In~\cref{sec: semantics for hierarchical nets} we
defined hierarchical nets as nets endowed with a specific kind of
functorial semantics to $\Span$. As things stand now, Petri nets
correspond to categories, while hierarchical nets correspond to
functors. This difference makes it difficult to say what a Petri net
with multiple levels of hierarchy is: intuitively, it is easy to imagine that
the children of a parent net $N$ can be themselves parents of other
nets, which are thus ``grandchildren'' of $N$, and so on and so forth.

In realizing this, we are blocked by having to map $N$ to hierarchical nets,
which are functors and not categories. To make such an intuition viable,
we need a way to \emph{internalize} the semantics in~\cref{def: hierarchical nets external definition}
to obtain a category representing the executions of the hierarchical net.

Luckily, there is a way to turn functors into categories, which relies on
an equivalence between the slice 2-category over a given category $\CategoryC$, denoted $\Cat/\CategoryC$, and the 2-category of lax-functors $\CategoryC \to \Span$~\cite{Pavlovic1997}.
This is itself the ``1-truncated'' version of a more general equivalence between
the slice of $\Cat$ over $\CategoryC$, and the 2-category of lax \emph{normal} functors to
the bicategory $\cate{Prof}$ of profunctors (this has been discovered by B\'enabou~\cite{Benabou1967}; a fully worked out exposition, conducted in full detail, is in~\cite{coend-calcu}).

\begin{figure}
  \begin{small}
    \begin{center}
        \scalebox{0.7}{
    \begin{tikzpicture}[node distance=1.3cm,>=stealth',bend angle=45,baseline=(current bounding box.center)]
      \node [place, tokens=0,] (1a) at (0,0){};
      \node [transition] (2a) at (1.5,0)     {}
          edge [pre] (1a);
      \node [place,tokens=0] (3a) at (3,0)  {}
          edge [pre] (2a);
      \node [transition] (4a) at (4.5,0) {}
          edge [pre] (3a);
      \node [place,tokens=0] (5a) at (6,0) {}
          edge [pre] (4a);
      \draw[thick, draw=black!75, fill=black!20] (0,2) ellipse (0.25cm and 0.75cm);
      \draw[thick, draw=black!75, fill=black!20] (3,2) ellipse (0.25cm and 0.75cm);
      \draw[thick, draw=black!75, fill=black!20] (6,2) ellipse (0.25cm and 0.75cm);

      \node[circle,fill=red, draw=gray!75, inner sep=0pt,minimum size=4pt] (red) at (0,1.7) {};
      \node[circle,fill=blue, draw=gray!75, inner sep=0pt,minimum size=4pt] (blue) at (0,2.3) {};

      \node[circle,fill=yellow, draw=gray!75, inner sep=0pt,minimum size=4pt] (yellow) at (3,1.7) {};
      \node[circle,fill=green, draw=gray!75, inner sep=0pt,minimum size=4pt] (green) at (3,2.3) {};

      \node[circle,fill=brown, draw=gray!75, inner sep=0pt,minimum size=4pt] (brown) at (5.9,2.4) {};
      \node[circle,fill=purple, draw=gray!75, inner sep=0pt,minimum size=4pt] (purple) at (6.1,2) {};
      \node[circle,fill=orange, draw=gray!75, inner sep=0pt,minimum size=4pt] (orange) at (5.9,1.6) {};

      \draw[->, out=20, in=160] (blue) to (green);
      \draw[->, out=0, in=200] (red) to (green);
      \draw[->, out=0, in=180] (yellow) to (purple);

      \draw[dotted, -]  (1a.north) -- (0,1.25);
      \draw[dotted, -]  (3a.north) -- (3,1.25);
      \draw[dotted, -]  (5a.north) -- (6,1.25);

      \node[transition,dotted,fill=none, draw=black!50, inner sep=0pt,minimum height=40pt] (left) at (1.5,2.2) {};
      \node[transition,dotted,fill=none, draw=black!50, inner sep=0pt,minimum height=10pt] (right) at (4.5,1.85) {};

      \draw[dotted, -]  (2a.north) -- (left.south);
      \draw[dotted, -]  (4a.north) -- (right.south);

    \end{tikzpicture}
  }
  $\qquad \leadsto \qquad$
  \scalebox{0.7}{
    \begin{tikzpicture}[node distance=1.3cm,>=stealth',bend angle=45, baseline=(current bounding box.center)]
      \node [place, tokens=0,] (1a) at (0,0){};
      \node [transition] (2a) at (1.5,0)     {}
          edge [pre] (1a);
      \node [place,tokens=0] (3a) at (3,0)  {}
          edge [pre] (2a);
      \node [transition] (4a) at (4.5,0) {}
          edge [pre] (3a);
      \node [place,tokens=0] (5a) at (6,0) {}
          edge [pre] (4a);
      \draw[thick, dotted, draw=black!50, fill=black!10] (0,2) ellipse (0.5cm and 1.25cm);
      \draw[thick, dotted, draw=black!50, fill=black!10] (3,2) ellipse (0.5cm and 1.25cm);
      \draw[thick, dotted, draw=black!50, fill=black!10] (6,2) ellipse (0.8cm and 1.25cm);

      \node[place, draw=red!50, tokens=0] (red) at (0,1.5) {};
      \node[place, draw=blue!50, tokens=0] (blue) at (0,2.5) {};

      \node[place, draw=yellow!50, tokens=0] (yellow) at (3,1.5) {};
      \node[place, draw=green!50, tokens=0] (green) at (3,2.5) {};

      \node[place, draw=brown!50, tokens=0] (brown) at (5.8,2.7) {};
      \node[place, draw=purple!50, tokens=0] (purple) at (6.2,2) {};
      \node[place, draw=orange!50, tokens=0] (orange) at (5.8,1.3) {};

      \draw[dotted, -]  (2a.north) -- (1.5,3);
      \draw[dotted, -]  (4a.north) -- (4.5,1.8);

      \draw[->, out=20, in=160] (blue) to node[transition] {} (green);
      \draw[->, out=0, in=200] (red) to node[transition] {} (green);
      \draw[->, out=0, in=180] (yellow) to node[pos=0.455, transition] {} (purple);

      \draw[dotted, -]  (1a.north) -- (0,0.75);
      \draw[dotted, -]  (3a.north) -- (3,0.75);
      \draw[dotted, -]  (5a.north) -- (6,0.75);

    \end{tikzpicture}
  }
    \end{center}
    \caption{The Grothendieck construction applied to $\Fun{N}$.}
    \label{fig:grote}
  \end{small}
\end{figure}
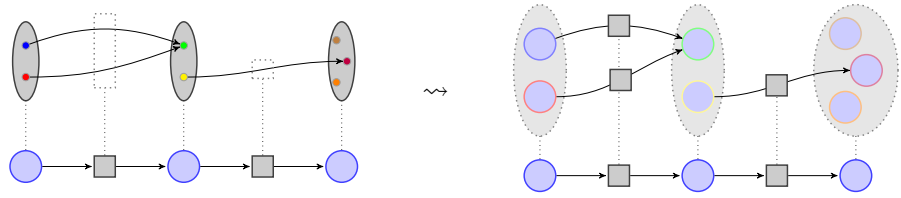

Here, we gloss over these abstract motivations and just give a very explicit
definition of what this means, as what we need is just a particular case of the construction
we worked out for guarded nets in~\cite{Genovese2020}.
\begin{definition}[Internalization]\label{def: internalization}
  Let $\NetSem{M}\in\PetriHier$ be a hierarchical net. We define
  its \emph{internalization}, denoted $\Grothendieck{\Fun{M}}$, as the following category:
  \begin{itemize}
    \item The objects of $\GrothendieckS{M}$ are pairs
    $(X, x)$, where $X$ is an object of $\Free{M}$ and $x$ is an
    element of $\Fun{M}X$. Concisely:
    \begin{equation*}
      \Obj{(\GrothendieckS{M})} :=
        \Suchthat{(X,x)}{(X \in \Obj{\Free{M}}) \wedge (x \in \Fun{M}X)}.
    \end{equation*}
    \item A morphism from $(X,x)$ to
    $(Y,y)$ in $\Grothendieck{\Fun{M}}$ is a pair $(f,s)$ where $f\colon X \to Y$ in $\Free{M}$ and
    $s\in S_{\Fun{M}f}$ in the apex of the corresponding span that connects $x$ to $y$. Concisely:
    \begin{align*}
      &\Hom{\GrothendieckS{M}}{(X,x)}{(Y,y)} :=\\
       &\qquad :=\Suchthat{(f,s)}{(f \in \Hom{\Free{M}}{X}{Y}) \wedge (s\in S_{\Fun{M}f})\wedge(\Fun{M}f(s) = (x,y))}.
    \end{align*}
  \end{itemize}
\end{definition}
The category $\GrothendieckS{N}$, called \emph{the Grothendieck construction applied
to $\Fun{N}$}, produces a place for each element of the set we send a place
to, and makes a transition for each path between these elements,
as shown in Figure \ref{fig:grote}.

Notice that in~\cref{fig:grote}, on the left, each path between coloured
dots is a triple $(a,x,b)$ as in~\cref{rem: interpreting hierarchical nets}.
This amounts to promote every possible trace of the child net -- together
with a selection of initial and accepting states -- to a transition in the parent
net. This interpretation is justified by the following theorem, which we again
proved in~\cite{Genovese2020}:
\begin{theorem}\label{thm: internalization}
  Given any strict monoidal functor $\Free{N} \xrightarrow{\Fun{N}} \Span$, the category $\GrothendieckS{N}$ is
  symmetric strict monoidal, and free. Thus $\GrothendieckS{N}$ can be written as $\Free{M}$ for some net $M$.

  Moreover, we obtain a \emph{projection functor} $\GrothendieckS{N} \to \Free{N}$ which turns $\Grothendieck{}$
  into a functor, in that for each functor $F: \NetSem{M} \to \NetSem{N}$ there exists a functor $\LiftSpan{F}$ making
  the following diagram commute:
  \begin{equation*}
    \begin{tikzpicture}[node distance=1.3cm,>=stealth',bend angle=45,auto]
      \node (1) at (0,1.5) {$\GrothendieckS{M}$};
      \node (2) at (0,0) {$\Free{M}$};
      \node (3) at (3,1.5) {$\GrothendieckS{N}$};
      \node (4) at (3,0) {$\Free{N}$};
      \node (5) at (1.5, -1.5) {$\Span$};
      \draw[->] (1)--(2) node [midway,left] {$\pi_M$};
      \draw[->] (3)--(4) node [midway,right] {$\pi_N$};
      \draw[->] (2)--(4) node [midway,above] {$F$};
      \draw[->, dashed] (1)--(3) node [midway,above] {$\LiftSpan{F}$};
      \draw[->] (2)--(5) node [midway,left] {$\Fun{M}$};
      \draw[->] (4)--(5) node [midway,right] {$\Fun{N}$};
    \end{tikzpicture}
  \end{equation*}
\end{theorem}
\cref{thm: internalization} defines a functor $\PetriSpan \to \FSSMC$,
the category of {\fssmc}s and strict monoidal functors between them. As $\PetriHier$
is a subcategory of $\PetriSpan$, we can immediately restrict~\cref{thm: internalization}
to hierarchical nets. A net in the form $\GrothendieckS{N}$ for some hiearchical net $\NetSem{N}$
is called the \emph{internal categorical semantics for $N$} (compare this
with~\cref{def: hierarchical nets external definition}, which we called \emph{external}).
\begin{remark}
  Notice how internalization is \emph{very} different from just copy-pasting a child net in place of a transition in the parent net as we discussed in~\cref{sec: hierarchical nets}.
  Here, each \emph{execution} of the child net is promoted to a transition, preserving the atomicity requirement of transitions in the parent net.
\end{remark}
Clearly, now we can define hierarchical nets with a level of hierarchy higher than two by just mapping a
generator $f$ of the parent net to a span where $N_f$ is in the form $\GrothendieckS{N}$ for some other hierarchical nets $N$, and the process can be recursively applied any finite number of times for each
transition.

\section{Engineering perspective}\label{sec: engineering perspective}
We deem it wise to spend a few words on why we consider this way of doing
things advantageous from an applicative perspective. Petri nets have been
considered as a possible way of producing software for a long time, with some
startups even using them as a central tool in their product offer~\cite{StateboxTeam2019}. Providing
some form of hierarchical calling is needed to make the idea of ``Petri
nets as a programming language/general-purpose design tool'' practical.

Our definition of hierarchy has the advantage of not making hierarchical
nets more expressive than Petri nets. If this seems like a downside, notice that
a consequence of this is that decidability of any reachability-related question is exactly
as for Petri nets, which is a great advantage from the point of view of model checking.
The legitimacy of this assertion is provided by internalization, that allows us to 
reduce hierarchical nets back to Petri nets. A further advantage of this is that
we can use already widespread tools for reachability checking~\cite{UniversityofTorino2018}
to answer reachability questions for our hierarchical nets,
without having necessarily to focus on producing new ones.

Moreover, and more importantly, our span formalism works really well
in modelling net behaviour in a distributed setting. 
To better understand this, imagine an infrastructure where each Petri net is
considered as a \emph{smart contract} (as it would be, for instance,
if we were to implement nets as smart contracts on a blockchain). A smart contract is nothing
more than a piece of code residing at a given address. Interaction with smart contracts
is \emph{transactional}: One sends a request to the contract address with some data to
be processed (for example a list of functions to be called on some parameters). The smart contract
executes as per the transaction and returns the data processed.

In our Petri net example things do not change: A user sends a
message consisting of a net address, the transaction the
user intends to fire, and some transaction data. The infrastructure replies affirmatively
or negatively if the transaction can be fired, which amounts to accept or reject the transaction.
As we already stressed, this is particularly suitable for blockchain-related
contexts and it is how applications such as~\cite{StateboxTeam2017a} implement Petri nets in their services.
\begin{figure}[ht!]\centering
  \input{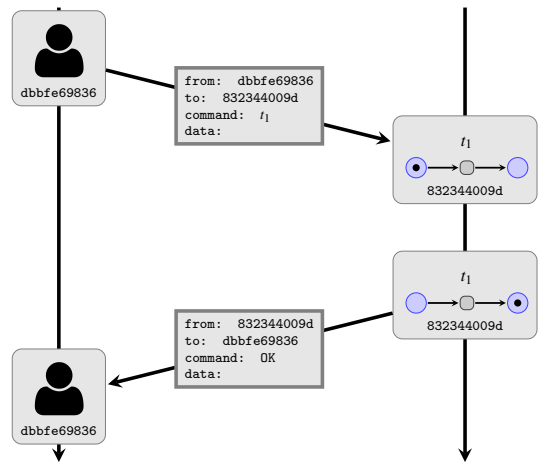}
  \caption{In this diagram we describe the interaction between a user and a net, with downward pointing arrows
  representing the flow of time. The user, having id \texttt{dbbfe69836}, sends a request to a net having address \texttt{832344009d}.
  The user is requesting to fire transition $t_1$ in the net. As the transition is enabled and able to fire, the request is granted,
  the state of the net updated, and a reply to the user is sent.}
\end{figure}

From this point of view, a hierarchical net would work exactly
as a standard Petri net, with the exception that in sending a transaction
to the parent net, the user also has to specify, in the transaction data,
a proper execution of the child net corresponding to the firing transition.
\begin{figure*}[ht!]\centering
  \input{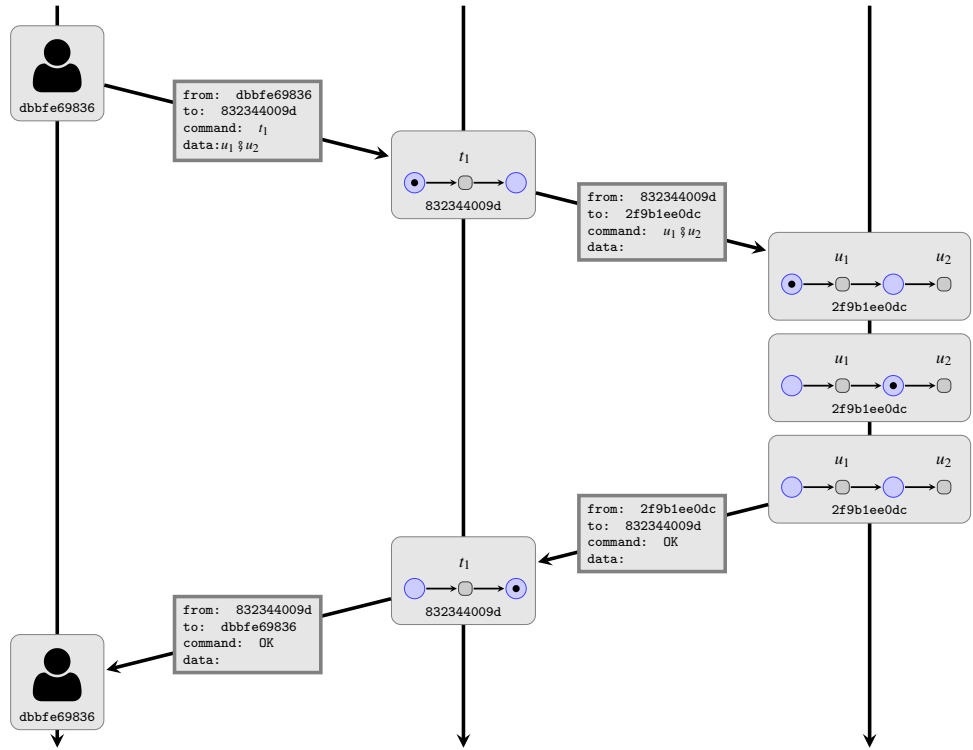} 
  \caption{In this diagram we describe the interaction between a user and a hierarchical net. 
  This time the user, having id \texttt{dbbfe69836}, sends a request to a net having address \texttt{832344009d}.
  This net is hierarchical, so in calling transition $t_1$ in the parent net, the user has also to provide a valid execution for
  its child. This is provided as transaction data, in this case $u_1 \cp u_2$. The parent net stores the address to the
  child net corresponding to $t_1$, which in this case is \texttt{2f9b1ee0dc}. The request to fire $u_1$ and then $u_2$ is
  forwarded to \texttt{2f9b1ee0dc}, which changes its state and responds affirmatively. This means that \texttt{832344009d}
  can itself change its state and respond affirmatively to \texttt{dbbfe69836}. Should any of these steps fail, the entire
  transaction is rejected and each net reverts to its previous state.}
\end{figure*}

Again, from a smart contract standpoint, this means that the smart contract
corresponding to the parent net will call the contract corresponding to the child
net with some execution data, and will respond affirmatively to the user only if
the generated call resolves positively.

Recalling the results in previous sections of this work, 
all the possible ways of executing the contracts above
form a category, which is obtained by internalizing the
hierarchical net via~\cref{thm: internalization}.
Internalized categories being free, they are presented by
Petri nets, which we can feed to any mainstream model
checker. Now, all sorts of questions about liveness
and interaction of the contracts above can be analyzed
by model-checking the corresponding internalized net.
This provides an easy way to analyze complex contract
interaction, relying on tools that have been debugged and
computationally optimized for decades.

\section{Conclusion and future work}\label{sec: conclusion and future work}
In this work, we showed how a formalism for guarded nets already worked out in~\cite{Genovese2020}
can be used to define the categorical semantics of some particular variety of hierarchical
nets, which works particularly well from a model-checking and distributed-implementation
point of view. Our effort is again part of a more ample project focusing on characterizing
the categorical semantics of extensions of Petri nets by studying functors from {\fssmc}s
to spans~\cite{Genovese2021a, Genovese2021}.

As a direction of future work, we would like to obtain a cleaner way of describing
recursively hierarchical nets. In this work, we relied on the Grothendieck construction
to internalize a hierarchical net, so that we could use hierarchical nets as children of
some other another parent net, recursively. This feels a bit like throwing all the
carefully-typed information that the external semantics gives into the same bucket,
and as such it is a bit unsatisfactory. Ideally, we would like to get a fully external
semantics for recursively hierarchical nets, and generalize the internalization result
to this case.

Another obvious direction of future work is implementing the findings hereby
presented, maybe relying on some formally verified implementation of category
theory such as~\cite{Genovese2019d}.
\paragraph*{\bf Acknowledgements}
Being non-native speakers, we want to thank John Baez for pointing out that the correct spelling in English was ``Hierarchical Petri nets" and not ``Hierarchic Petri nets" as we thought. We also want to thank G. Pasquini and M. De Jorio for having inspired this work.

\bigskip\noindent
A video presentation of this paper can be found on Youtube at \href{https://www.youtube.com/watch?v=4v5v8tgmiUM}{4v5v8tgmiUM}.

\bibliographystyle{eptcs}
  \bibliography{Bibliography}

\begin{thebibliography}{10}
\providecommand{\bibitemdeclare}[2]{}
\providecommand{\surnamestart}{}
\providecommand{\surnameend}{}
\providecommand{\urlprefix}{Available at }
\providecommand{\url}[1]{\texttt{#1}}
\providecommand{\href}[2]{\texttt{#2}}
\providecommand{\urlalt}[2]{\href{#1}{#2}}
\providecommand{\doi}[1]{doi:\urlalt{http://dx.doi.org/#1}{#1}}
\providecommand{\bibinfo}[2]{#2}

\bibitemdeclare{inproceedings}{Baez2021}
\bibitem{Baez2021}
\bibinfo{author}{John~C. \surnamestart Baez\surnameend},
  \bibinfo{author}{Fabrizio \surnamestart Genovese\surnameend},
  \bibinfo{author}{Jade \surnamestart Master\surnameend} \&
  \bibinfo{author}{Michael \surnamestart Shulman\surnameend}
  (\bibinfo{year}{2021}): \emph{\bibinfo{title}{Categories of Nets}}.
\newblock In: {\sl \bibinfo{booktitle}{36th Annual {ACM/IEEE} Symposium on
  Logic in Computer Science, {LICS} 2021, Rome, Italy, June 29 - July 2,
  2021}}, \bibinfo{publisher}{{IEEE}}, pp. \bibinfo{pages}{1--13},
  \doi{10.1109/LICS52264.2021.9470566}.

\bibitemdeclare{article}{Baez2020a}
\bibitem{Baez2020a}
\bibinfo{author}{John~C. \surnamestart Baez\surnameend} \&
  \bibinfo{author}{Jade \surnamestart Master\surnameend}
  (\bibinfo{year}{2020}): \emph{\bibinfo{title}{Open Petri nets}}.
\newblock {\sl \bibinfo{journal}{Mathematical Structures in Computer Science}}
  \bibinfo{volume}{30}(\bibinfo{number}{3}), p. \bibinfo{pages}{314–341},
  \doi{10.1017/S0960129520000043}.

\bibitemdeclare{article}{baldan2018petri}
\bibitem{baldan2018petri}
\bibinfo{author}{P.~\surnamestart Baldan\surnameend},
  \bibinfo{author}{M.~\surnamestart Bocci\surnameend},
  \bibinfo{author}{D.~\surnamestart Brigolin\surnameend},
  \bibinfo{author}{N.~\surnamestart Cocco\surnameend},
  \bibinfo{author}{M.~\surnamestart Heiner\surnameend} \&
  \bibinfo{author}{M.~\surnamestart Simeoni\surnameend} (\bibinfo{year}{2018}):
  \emph{\bibinfo{title}{Petri nets for modelling and analysing trophic
  networks}}.
\newblock {\sl \bibinfo{journal}{Fundamenta Informaticae}}
  \bibinfo{volume}{160}(\bibinfo{number}{1-2}), pp. \bibinfo{pages}{27--52},
  \doi{10.3233/FI-2018-1673}.

\bibitemdeclare{inproceedings}{baldan2014encoding}
\bibitem{baldan2014encoding}
\bibinfo{author}{P.~\surnamestart Baldan\surnameend},
  \bibinfo{author}{F.~\surnamestart Bonchi\surnameend},
  \bibinfo{author}{F.~\surnamestart Gadducci\surnameend} \&
  \bibinfo{author}{G.V. \surnamestart Monreale\surnameend}
  (\bibinfo{year}{2014}): \emph{\bibinfo{title}{Encoding synchronous
  interactions using labelled Petri nets}}.
\newblock In: {\sl \bibinfo{booktitle}{International Conference on Coordination
  Languages and Models}}, \bibinfo{organization}{Springer}, pp.
  \bibinfo{pages}{1--16}, \doi{10.1007/978-3-662-43376-8\_1}.

\bibitemdeclare{incollection}{baldan2015asynchronous}
\bibitem{baldan2015asynchronous}
\bibinfo{author}{P.~\surnamestart Baldan\surnameend},
  \bibinfo{author}{F.~\surnamestart Bonchi\surnameend},
  \bibinfo{author}{F.~\surnamestart Gadducci\surnameend} \&
  \bibinfo{author}{G.V. \surnamestart Monreale\surnameend}
  (\bibinfo{year}{2015}): \emph{\bibinfo{title}{Asynchronous Traces and Open
  Petri Nets}}.
\newblock In: {\sl \bibinfo{booktitle}{Programming Languages with Applications
  to Biology and Security}}, \bibinfo{publisher}{Springer}, pp.
  \bibinfo{pages}{86--102}, \doi{10.1007/978-3-319-25527-9\_8}.

\bibitemdeclare{article}{baldan2015modular}
\bibitem{baldan2015modular}
\bibinfo{author}{P.~\surnamestart Baldan\surnameend},
  \bibinfo{author}{F.~\surnamestart Bonchi\surnameend},
  \bibinfo{author}{F.~\surnamestart Gadducci\surnameend} \&
  \bibinfo{author}{G.V. \surnamestart Monreale\surnameend}
  (\bibinfo{year}{2015}): \emph{\bibinfo{title}{Modular encoding of synchronous
  and asynchronous interactions using open Petri nets}}.
\newblock {\sl \bibinfo{journal}{Science of Computer Programming}}
  \bibinfo{volume}{109}, pp. \bibinfo{pages}{96--124},
  \doi{10.1016/j.scico.2014.11.019}.

\bibitemdeclare{inproceedings}{baldan2011mpath2pn}
\bibitem{baldan2011mpath2pn}
\bibinfo{author}{P.~\surnamestart Baldan\surnameend},
  \bibinfo{author}{N.~\surnamestart Cocco\surnameend},
  \bibinfo{author}{F.~\surnamestart De~Nes\surnameend}, \bibinfo{author}{M.L.
  \surnamestart Segura\surnameend} \& \bibinfo{author}{M.~\surnamestart
  Simeoni\surnameend} (\bibinfo{year}{2011}):
  \emph{\bibinfo{title}{MPath2PN-Translating metabolic pathways into Petri
  nets}}.
\newblock In: {\sl \bibinfo{booktitle}{BioPPN2011 Int. Workshop on Biological
  Processes and Petri Nets, CEUR Workshop Proceedings}}, \bibinfo{volume}{724},
  pp. \bibinfo{pages}{102--116}.

\bibitemdeclare{inproceedings}{baldan2001compositional}
\bibitem{baldan2001compositional}
\bibinfo{author}{P.~\surnamestart Baldan\surnameend},
  \bibinfo{author}{A.~\surnamestart Corradini\surnameend},
  \bibinfo{author}{H.~\surnamestart Ehrig\surnameend} \&
  \bibinfo{author}{R.~\surnamestart Heckel\surnameend} (\bibinfo{year}{2001}):
  \emph{\bibinfo{title}{Compositional modeling of reactive systems using open
  nets}}.
\newblock In: {\sl \bibinfo{booktitle}{International Conference on Concurrency
  Theory}}, \bibinfo{organization}{Springer}, pp. \bibinfo{pages}{502--518},
  \doi{10.1007/3-540-44685-0\_34}.

\bibitemdeclare{inproceedings}{baldan2008open}
\bibitem{baldan2008open}
\bibinfo{author}{P.~\surnamestart Baldan\surnameend},
  \bibinfo{author}{A.~\surnamestart Corradini\surnameend},
  \bibinfo{author}{H.~\surnamestart Ehrig\surnameend} \&
  \bibinfo{author}{B.~\surnamestart K{\"o}nig\surnameend}
  (\bibinfo{year}{2008}): \emph{\bibinfo{title}{Open Petri nets:
  Non-deterministic processes and compositionality}}.
\newblock In: {\sl \bibinfo{booktitle}{International Conference on Graph
  Transformation}}, \bibinfo{organization}{Springer}, pp.
  \bibinfo{pages}{257--273}, \doi{10.1007/978-3-540-87405-8\_18}.

\bibitemdeclare{article}{baldan2010petri}
\bibitem{baldan2010petri}
\bibinfo{author}{P.~\surnamestart Baldan\surnameend},
  \bibinfo{author}{A.~\surnamestart Corradini\surnameend},
  \bibinfo{author}{F.~\surnamestart Gadducci\surnameend} \&
  \bibinfo{author}{U.~\surnamestart Montanari\surnameend}
  (\bibinfo{year}{2010}): \emph{\bibinfo{title}{From Petri nets to graph
  transformation systems}}.
\newblock {\sl \bibinfo{journal}{Electronic Communications of the EASST}}
  \bibinfo{volume}{26}, \doi{10.14279/tuj.eceasst.26.368}.

\bibitemdeclare{article}{baldan2005relating}
\bibitem{baldan2005relating}
\bibinfo{author}{P.~\surnamestart Baldan\surnameend},
  \bibinfo{author}{A.~\surnamestart Corradini\surnameend} \&
  \bibinfo{author}{U.~\surnamestart Montanari\surnameend}
  (\bibinfo{year}{2005}): \emph{\bibinfo{title}{Relating SPO and DPO graph
  rewriting with Petri nets having read, inhibitor and reset arcs}}.
\newblock {\sl \bibinfo{journal}{Electronic Notes in Theoretical Computer
  Science}} \bibinfo{volume}{127}(\bibinfo{number}{2}), pp.
  \bibinfo{pages}{5--28}, \doi{10.1016/j.entcs.2005.02.003}.

\bibitemdeclare{article}{baldan2019petri}
\bibitem{baldan2019petri}
\bibinfo{author}{P.~\surnamestart Baldan\surnameend} \&
  \bibinfo{author}{F.~\surnamestart Gadducci\surnameend}
  (\bibinfo{year}{2019}): \emph{\bibinfo{title}{Petri nets are dioids: a new
  algebraic foundation for non-deterministic net theory}}.
\newblock {\sl \bibinfo{journal}{Acta Informatica}}
  \bibinfo{volume}{56}(\bibinfo{number}{1}), pp. \bibinfo{pages}{61--92},
  \doi{10.1007/s00236-018-0314-0}.

\bibitemdeclare{inbook}{Benabou1967}
\bibitem{Benabou1967}
\bibinfo{author}{J.~\surnamestart B\'enabou\surnameend}:
  \emph{\bibinfo{title}{Introduction to {{Bicategories}}}}, pp.
  \bibinfo{pages}{1--77}.
\newblock \bibinfo{volume}{47}, \bibinfo{publisher}{{Springer Berlin
  Heidelberg}}, \doi{10.1007/BFb0074299}.

\bibitemdeclare{inproceedings}{Buchholz1994HierarchicalHL}
\bibitem{Buchholz1994HierarchicalHL}
\bibinfo{author}{P.~\surnamestart Buchholz\surnameend} (\bibinfo{year}{1994}):
  \emph{\bibinfo{title}{Hierarchical High Level Petri Nets for Complex System
  Analysis}}.
\newblock In: {\sl \bibinfo{booktitle}{Application and Theory of Petri Nets}},
  \doi{10.1007/3-540-58152-9\_8}.

\bibitemdeclare{unpublished}{benabou2000distributors}
\bibitem{benabou2000distributors}
\bibinfo{author}{J.~\surnamestart Bénabou\surnameend} \&
  \bibinfo{author}{T.~\surnamestart Streicher\surnameend}
  (\bibinfo{year}{2000}): \emph{\bibinfo{title}{Distributors at work}}.
\newblock \bibinfo{note}{Lecture notes written by Thomas Streicher}.

\bibitemdeclare{article}{cattani2005profunctors}
\bibitem{cattani2005profunctors}
\bibinfo{author}{G.~L. \surnamestart Cattani\surnameend} \&
  \bibinfo{author}{G.~\surnamestart Winskel\surnameend} (\bibinfo{year}{2005}):
  \emph{\bibinfo{title}{Profunctors, open maps and bisimulation}}.
\newblock {\sl \bibinfo{journal}{Mathematical Structures in Computer Science}}
  \bibinfo{volume}{15}(\bibinfo{number}{03}), pp. \bibinfo{pages}{553--614},
  \doi{10.1017/S0960129505004718}.

\bibitemdeclare{article}{Esparza1994}
\bibitem{Esparza1994}
\bibinfo{author}{J.~\surnamestart Esparza\surnameend} \&
  \bibinfo{author}{N.~\surnamestart Mogens\surnameend} (\bibinfo{year}{1994}):
  \emph{\bibinfo{title}{Decidability Issues for Petri Nets - a survey}}.
\newblock {\sl \bibinfo{journal}{J. Inf. Process. Cybern.}}
  \bibinfo{volume}{30}(\bibinfo{number}{3}), pp. \bibinfo{pages}{143--160}.

\bibitemdeclare{inproceedings}{fehling1991concept}
\bibitem{fehling1991concept}
\bibinfo{author}{R.~\surnamestart Fehling\surnameend} (\bibinfo{year}{1991}):
  \emph{\bibinfo{title}{A concept of hierarchical Petri nets with building
  blocks}}.
\newblock In: {\sl \bibinfo{booktitle}{International Conference on Application
  and Theory of Petri Nets}}, \bibinfo{organization}{Springer}, pp.
  \bibinfo{pages}{148--168}, \doi{10.1007/3-540-56689-9\_43}.

\bibitemdeclare{}{Genovese2019d}
\bibitem{Genovese2019d}
\bibinfo{author}{F.~\surnamestart Genovese\surnameend},
  \bibinfo{author}{A.~\surnamestart Gryzlov\surnameend},
  \bibinfo{author}{J.~\surnamestart Herold\surnameend},
  \bibinfo{author}{A.~\surnamestart Knispel\surnameend},
  \bibinfo{author}{M.~\surnamestart Perone\surnameend},
  \bibinfo{author}{E.~\surnamestart Post\surnameend} \&
  \bibinfo{author}{A.~\surnamestart Videla\surnameend}:
  \emph{\bibinfo{title}{Idris-Ct: {{A Library}} to Do {{Category Theory}} in
  {{Idris}}}}.
\newblock \urlprefix\url{http://arxiv.org/abs/1912.06191}.

\bibitemdeclare{}{Genovese2019c}
\bibitem{Genovese2019c}
\bibinfo{author}{F.~\surnamestart Genovese\surnameend},
  \bibinfo{author}{A.~\surnamestart Gryzlov\surnameend},
  \bibinfo{author}{J.~\surnamestart Herold\surnameend},
  \bibinfo{author}{M.~\surnamestart Perone\surnameend},
  \bibinfo{author}{E.~\surnamestart Post\surnameend} \&
  \bibinfo{author}{A.~\surnamestart Videla\surnameend}:
  \emph{\bibinfo{title}{Computational {{Petri Nets}}: {{Adjunctions Considered
  Harmful}}}}.
\newblock \urlprefix\url{http://arxiv.org/abs/1904.12974}.

\bibitemdeclare{article}{Genovese2019b}
\bibitem{Genovese2019b}
\bibinfo{author}{F.~\surnamestart Genovese\surnameend} \&
  \bibinfo{author}{J.~\surnamestart Herold\surnameend}:
  \emph{\bibinfo{title}{Executions in ({{Semi}}-){{Integer Petri Nets}} Are
  {{Compact Closed Categories}}}} \bibinfo{volume}{287}, pp.
  \bibinfo{pages}{127--144}.
\newblock \doi{10.4204/EPTCS.287.7}.

\bibitemdeclare{}{Genovese2021}
\bibitem{Genovese2021}
\bibinfo{author}{F.~\surnamestart Genovese\surnameend},
  \bibinfo{author}{F.~\surnamestart Loregian\surnameend} \&
  \bibinfo{author}{D.~\surnamestart Palombi\surnameend}:
  \emph{\bibinfo{title}{A {{Categorical Semantics}} for {{Bounded Petri
  Nets}}}}.
\newblock \urlprefix\url{http://arxiv.org/abs/2101.09100}.

\bibitemdeclare{}{Genovese2021a}
\bibitem{Genovese2021a}
\bibinfo{author}{F.~\surnamestart Genovese\surnameend},
  \bibinfo{author}{F.~\surnamestart Loregian\surnameend} \&
  \bibinfo{author}{D.~\surnamestart Palombi\surnameend}:
  \emph{\bibinfo{title}{Nets with {{Mana}}: {{A Framework}} for {{Chemical
  Reaction Modelling}}}}, \doi{10.1007/978-3-030-78946-6\_10}.
\newblock \urlprefix\url{http://arxiv.org/abs/2101.06234}.

\bibitemdeclare{incollection}{Genovese2020}
\bibitem{Genovese2020}
\bibinfo{author}{F.~\surnamestart Genovese\surnameend} \& \bibinfo{author}{D.I.
  \surnamestart Spivak\surnameend}: \emph{\bibinfo{title}{A {{Categorical
  Semantics}} for {{Guarded Petri Nets}}}}.
\newblock In \bibinfo{editor}{F.~\surnamestart Gadducci\surnameend} \&
  \bibinfo{editor}{Timo \surnamestart Kehrer\surnameend}, editors: {\sl
  \bibinfo{booktitle}{Graph {{Transformation}}}}, {\sl \bibinfo{series}{Lecture
  {{Notes}} in {{Computer Science}}}} \bibinfo{volume}{12150},
  \bibinfo{publisher}{{Springer International Publishing}}, pp.
  \bibinfo{pages}{57--74}, \doi{10.1007/978-3-030-51372-6\_4}.

\bibitemdeclare{inproceedings}{huber1989hierarchies}
\bibitem{huber1989hierarchies}
\bibinfo{author}{P.~\surnamestart Huber\surnameend},
  \bibinfo{author}{K.~\surnamestart Jensen\surnameend} \& \bibinfo{author}{R.M.
  \surnamestart Shapiro\surnameend} (\bibinfo{year}{1989}):
  \emph{\bibinfo{title}{Hierarchies in coloured Petri nets}}.
\newblock In: {\sl \bibinfo{booktitle}{International Conference on Application
  and Theory of Petri Nets}}, \bibinfo{organization}{Springer}, pp.
  \bibinfo{pages}{313--341}, \doi{10.1007/978-3-662-06289-0\_3}.

\bibitemdeclare{book}{Jensen2009}
\bibitem{Jensen2009}
\bibinfo{author}{K.~\surnamestart Jensen\surnameend} \& \bibinfo{author}{L.M.
  \surnamestart Kristensen\surnameend}: \emph{\bibinfo{title}{Coloured {{Petri
  Nets}}}}.
\newblock \bibinfo{publisher}{{Springer Berlin Heidelberg}},
  \doi{10.1007/BFb0046842}.

\bibitemdeclare{article}{2catlimits}
\bibitem{2catlimits}
\bibinfo{author}{G.M. \surnamestart Kelly\surnameend} (\bibinfo{year}{1989}):
  \emph{\bibinfo{title}{Elementary observations on 2-categorical limits}}.
\newblock {\sl \bibinfo{journal}{Bulletin of the Australian Mathematical
  Society}} \bibinfo{volume}{39}, pp. \bibinfo{pages}{301--317},
  \doi{10.1017/S0004972700002781}.

\bibitemdeclare{article}{Kohler-Bussmeier2014}
\bibitem{Kohler-Bussmeier2014}
\bibinfo{author}{M.~\surnamestart K\"ohler-Bu{\ss}meier\surnameend}:
  \emph{\bibinfo{title}{A {{Survey}} of {{Decidability Results}} for
  {{Elementary Object Systems}}}} (\bibinfo{number}{1}), pp.
  \bibinfo{pages}{99--123}.
\newblock \doi{10.3233/FI-2014-983}.

\bibitemdeclare{book}{coend-calcu}
\bibitem{coend-calcu}
\bibinfo{author}{F.~\surnamestart Loregian\surnameend} (\bibinfo{year}{2021}):
  \emph{\bibinfo{title}{Coend Calculus}}.
\newblock {\sl \bibinfo{series}{London Mathematical Society Lecture Note
  Series}} \bibinfo{volume}{468}, \bibinfo{publisher}{Cambridge University
  Press}.
\newblock \bibinfo{note}{ISBN 9781108746120}.

\bibitemdeclare{article}{Master2020}
\bibitem{Master2020}
\bibinfo{author}{J.~\surnamestart Master\surnameend}:
  \emph{\bibinfo{title}{Petri {{Nets Based}} on {{Lawvere Theories}}}}
  \bibinfo{volume}{30}(\bibinfo{number}{7}), pp. \bibinfo{pages}{833--864}.
\newblock \doi{10.1017/S0960129520000262}.

\bibitemdeclare{article}{Meseguer1990}
\bibitem{Meseguer1990}
\bibinfo{author}{J.~\surnamestart Meseguer\surnameend} \&
  \bibinfo{author}{U.~\surnamestart Montanari\surnameend}:
  \emph{\bibinfo{title}{Petri {{Nets}} Are {{Monoids}}}}
  \bibinfo{volume}{88}(\bibinfo{number}{2}), pp. \bibinfo{pages}{105--155}.
\newblock \doi{10.1016/0890-5401(90)90013-8}.

\bibitemdeclare{inproceedings}{oswald1990environment}
\bibitem{oswald1990environment}
\bibinfo{author}{H.~\surnamestart Oswald\surnameend},
  \bibinfo{author}{R.~\surnamestart Esser\surnameend} \&
  \bibinfo{author}{R.~\surnamestart Mattmann\surnameend}
  (\bibinfo{year}{1990}): \emph{\bibinfo{title}{An environment for specifying
  and executing hierarchical Petri nets}}.
\newblock In: {\sl \bibinfo{booktitle}{[1990] Proceedings. 12th International
  Conference on Software Engineering}}, \bibinfo{organization}{IEEE}, pp.
  \bibinfo{pages}{164--172}, \doi{10.5555/100296.100319}.

\bibitemdeclare{incollection}{Pavlovic1997}
\bibitem{Pavlovic1997}
\bibinfo{author}{D.~\surnamestart Pavlovi\'c\surnameend} \&
  \bibinfo{author}{S.~\surnamestart Abramsky\surnameend}:
  \emph{\bibinfo{title}{Specifying {{Interaction Categories}}}}.
\newblock In \bibinfo{editor}{Eugenio \surnamestart Moggi\surnameend} \&
  \bibinfo{editor}{Giuseppe \surnamestart Rosolini\surnameend}, editors: {\sl
  \bibinfo{booktitle}{Category {{Theory}} and {{Computer Science}}}}, {\sl
  \bibinfo{series}{Lecture {{Notes}} in {{Computer Science}}}}
  \bibinfo{volume}{1290}, \bibinfo{publisher}{{Springer Berlin Heidelberg}},
  pp. \bibinfo{pages}{147--158}, \doi{10.5555/648335.755738}.

\bibitemdeclare{incollection}{Sassone1995}
\bibitem{Sassone1995}
\bibinfo{author}{V.~\surnamestart Sassone\surnameend}: \emph{\bibinfo{title}{On
  the {{Category}} of {{Petri Net Computations}}}}.
\newblock In: {\sl \bibinfo{booktitle}{{{TAPSOFT}} '95: {{Theory}} and
  {{Practice}} of {{Software Development}}}}, \bibinfo{volume}{915},
  \bibinfo{publisher}{{Springer Berlin Heidelberg}}, pp.
  \bibinfo{pages}{334--348}, \doi{10.1007/3-540-59293-8\_205}.

\bibitemdeclare{}{StateboxTeam2019}
\bibitem{StateboxTeam2019}
\bibinfo{author}{\surnamestart {Statebox Team}\surnameend}:
  \emph{\bibinfo{title}{The {{Mathematical Specification}} of the {{Statebox
  Language}}}}.
\newblock \urlprefix\url{http://arxiv.org/abs/1906.07629}.

\bibitemdeclare{}{StateboxTeam2017a}
\bibitem{StateboxTeam2017a}
\bibinfo{author}{\surnamestart {Statebox Team}\surnameend}:
  \emph{\bibinfo{title}{Statebox, {{Compositional Diagrammatic Programming
  Language}}}}.
\newblock \urlprefix\url{https://statebox.org}.

\bibitemdeclare{}{UniversityofTorino2018}
\bibitem{UniversityofTorino2018}
\bibinfo{author}{\surnamestart {University of Torino}\surnameend}:
  \emph{\bibinfo{title}{{{GreatSPN Github}} Page}}.
\newblock \urlprefix\url{https://github.com/greatspn/SOURCES}.

\bibitemdeclare{}{Zanasi2018}
\bibitem{Zanasi2018}
\bibinfo{author}{F.~\surnamestart Zanasi\surnameend}:
  \emph{\bibinfo{title}{Interacting {{Hopf Algebras}}: The {{Theory}} of
  {{Linear Systems}}}}.
\newblock \urlprefix\url{http://arxiv.org/abs/1805.03032}.

\end{thebibliography}
\end{document}